\documentclass[11pt]{amsart}
\usepackage{amsfonts,amssymb,amsbsy,amsmath,amscd,amstext,amsthm}

\usepackage[all]{xy}
\usepackage{amsmath}
\usepackage{amsthm}
\usepackage{amssymb}
\usepackage{amscd}
\usepackage{subfigure}
\usepackage{fancybox}
\usepackage{graphicx}
\usepackage{pstricks,graphicx}
\usepackage{epsfig}

\newtheorem{theorem}{Theorem}
\newtheorem{thm}[theorem]{Theorem}
\newtheorem{lem}[theorem]{Lemma}

\newtheorem{cor}[theorem]{Corollary}

\newtheorem{conjecture}{Theorem}
\newtheorem{conj}[conjecture]{Conjecture}

\theoremstyle{remark}

\def\N{\mathbb N}
\def\Z{\mathbb Z}
\def\Q{\mathbb Q}
\def\R{\mathbb R}

\def\C{\mathbb C}

\begin{document}

\title{Pentagonal Domain Exchange}
\author{Shigeki Akiyama}
\address{Department of Mathematics, Faculty of Science, Niigata University, Ikarashi-2 8050
Niigata, 950-2181 Japan}
\email{akiyama@math.sc.niigata-u.ac.jp}
\author{Edmund Harriss}
\address{Department of Mathematical Sciences, 1 University of Arkansas, 
Fayetteville, AR 72701, USA}
\email{edmund.harriss@mathematicians.org.uk}
\date{}
\thanks{The first author is supported by the Japanese Society for the Promotion of Science (JSPS), grant in aid 21540010.}
\maketitle

\begin{abstract}
Self-inducing structure of pentagonal piecewise isometry is applied to show detailed
description of periodic and aperiodic orbits, and further dynamical properties. 
A Pisot number appears as a scaling constant and plays a crucial role in the proof. 
Further generalization is discussed in the last section.
\end{abstract}

\begin{figure}[htbp]
	\centering
		\includegraphics{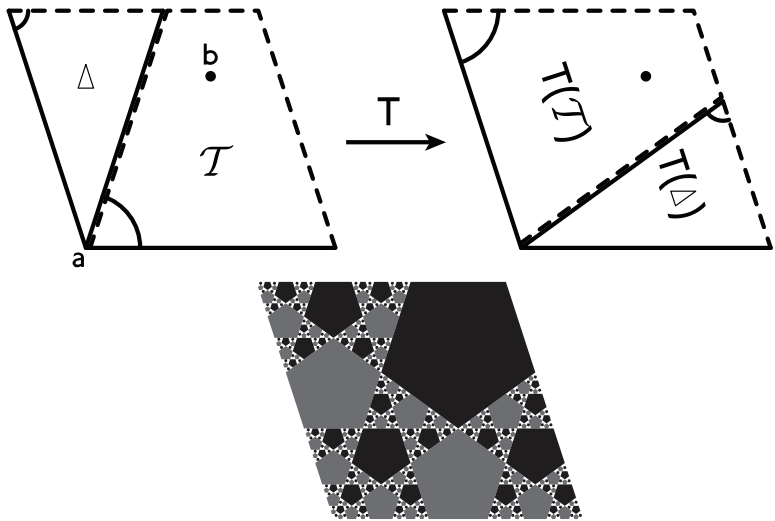}
	\caption{A piecewise rotation $T$ on two pieces. The triangle is rotated $2\pi/5$ around $a$ and the trapezium is rotated $2\pi/5$ around $b$. Periodic points with short periods are shown below, in two colours to illustrate that they cluster into groups, each forming a pentagon.}
	\label{Piece}
\end{figure}

Adler-Kitchens-Tresser~\cite{Adler-Kitchens-Tresser:01} and 
Goetz~\cite{Goetz:SEPIA} initiated the study of piecewise isometries.
This class of maps shows the way to possible 
generalizations of results on interval exchanges to higher 
dimensions~\cite{Kouptsov-Lowenstein-Vivaldi:02,Trovati:TPPTGPI}. 
In this paper we examine the detailed properties of the map 
shown in Figure \ref{Piece} from an algebraic point of view.

The goal of this paper is to see
how this map is applied to show number theoretical results.
First we reprove that almost all orbits 
in the sense of Lebesgue measure are periodic, and 
in addition, there are explicit aperiodic points. 
Second we show that
aperiodic points forms a proper dense subset 
of an attractor of some iterated function system and are
recognized by a B\"uchi automaton (c.f. Figure \ref{OpenEdge2}).
The dynamics acting on this set of aperiodic points
are conjugate to the $2$-adic odometer (addition of one)
whose explicit construction is given (Theorem \ref{Main}).
As a result, we easily see that all aperiodic orbits are dense and uniformly distributed 
in the attractor. We finally 
give a characterization of points which have purely periodic multiplicative coding
by constructing its natural extension (Theorem \ref{PureExpansion}).
In doing so we obtain an intriguing picture Figure \ref{Dual} that emerges naturally from taking algebraic conjugates, whose structure is worthy of further study. We
discuss possible generalizations for 7-fold and 9-fold piecewise rotations 
in Section \ref{sec:other_self_similar_systems}.

\begin{figure}[htbp]
	\centering
		\includegraphics{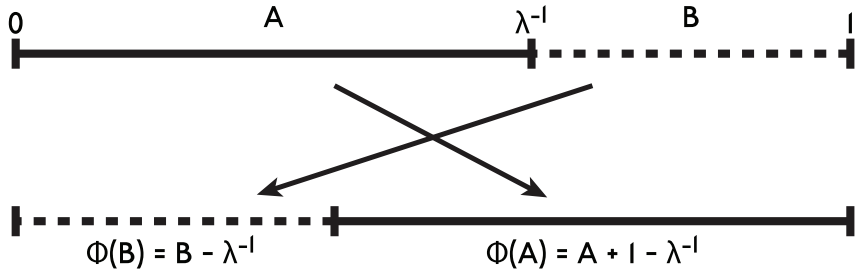}
	\caption{An interval exchange map $\Phi$, where $\lambda = \frac{1+\sqrt{5}}{2}$}
	\label{fig:Int_exchange}
\end{figure}

\begin{figure}[htbp]
	\centering
		\includegraphics{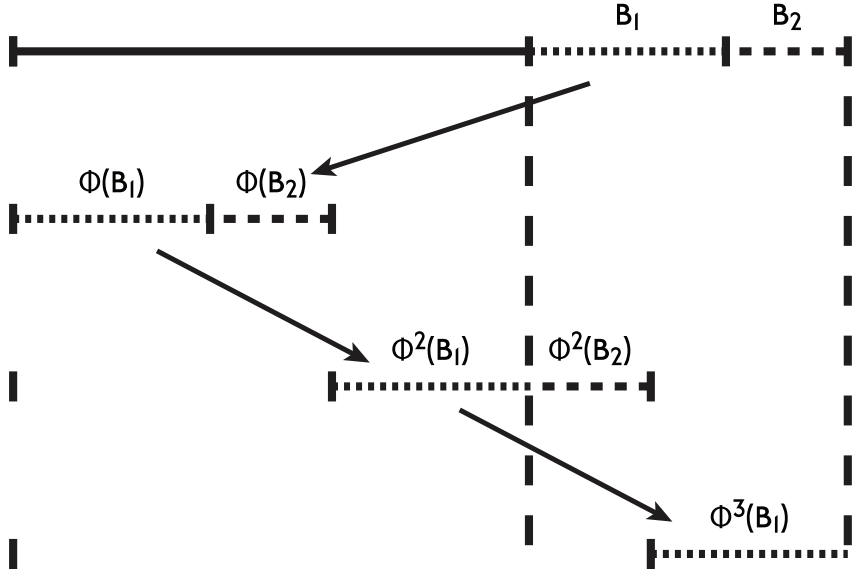}
	\caption{The interval exchange $\Phi$ is self-inducing. The intervals $B_1$ and $B_2$ are swapped by the first return map of $\Phi$ on the interval $B$.}
	\label{fig:Self_inducing}
\end{figure}

A dynamical system is {\it self-inducing} if the first return map to some subset has the same dynamics as the full map. The most important example is the irrational rotation, presented
as exchange of two intervals. An elementary example begins with $\Phi$ 
shown in Figure \ref{fig:Int_exchange}.
\noindent
For this interval exchange, now consider the second interval $B$. As shown in Figure \ref{fig:Self_inducing} this interval is translated to the left once, and to the right. Thus $\Phi^2(B_1)$ is back in $B$, the interval $B_2$ requires one more step, but $\Phi^3(B_2)$ also lies within $B$. This first return dynamics on $B$ is therefore conjugate to the dynamics on $A \cup B$.
\noindent
Self-inducing subsystem of two interval exchange 
corresponds to purely periodic orbits of 
continued fraction expansion
and they are efficiently captured by the continued fraction algorithm. 

This gives a motivation to study the interval exchange transform (IET) of three or
more pieces, trying to find higher dimensional continued fraction with good 
Diophantine approximation properties.
The study of self-inducing structure of IET's 
was started by a pioneer work of Rauzy \cite{Rauzy:EDTI}, now called Rauzy induction,
and got extended in
a great deal by many authors including Veech \cite{Veech:82}
and Zorich \cite{Zorich:96}, see \cite{Yoccoz:IET} for historical developments.

Self-inducing piecewise isometries emerged from dynamical systems as a natural generalization of IET~\cite{Lowenstein-Hatjispyros-Vivaldi,Adler-Kitchens-Tresser:01,Ashwin:GOPPI,Goetz:PIEADS,Goetz:RMCPS,Bressaud:TCBPPI,Mendes:SPPPIES}
and the first return dynamics appears in outer billiards \cite{Tabachnikov:95,
Bedaride-Cassaigne:10}.
Like IET they provide a simple setting to study many of the deep and perplexing behaviors that can emerge from a dynamical system. 

The self-inducing structure links such dynamical systems to 
number theoretical algorithms, such as, 
digital expansions and Diophantine approximation algorithms, and allows
us to study their periodic orbits by constructing their natural extensions.
This idea leads to complex and beautiful fractal behavior. 

Our target is the piecewise isometry 
in Figure \ref{Piece}, but to illustrate the bridge formed between the two 
fields let us begin with a simple 
conjecture from number theory:

\begin{conj}
\label{Rotation}
For any $-2<\lambda<2$, each integer sequence defined by
 $0\leq a_{n+1} + \lambda a_n + a_{n-1}<1$
is periodic.
\end{conj}

Since $a_{n+2}\in \Z$ is uniquely determined by $(a_n,a_{n+1})\in \Z^2$, we 
treat this recurrence as a map
$(a_n,a_{n+1})\mapsto (a_{n+1},a_{n+2})$ acting on $\Z^2$. 
It is natural to set $\lambda=-2\cos(\theta)$ to view this map 
as a `discretized rotation':
$$
\begin{pmatrix} a_{n+1}\\a_{n+2} 
\end{pmatrix} \sim \begin{pmatrix} 0 & 1\\ -1 & -\lambda\end{pmatrix}
\begin{pmatrix} a_n\\a_{n+1}\end{pmatrix}
$$
with eigenvalues $\exp(\pm \sqrt{-1} \theta)$. As the matrix is conjugate to 
the planar rotation matrix of angle $\theta$, putting $P=\begin{pmatrix}
1&0\cr \cos\theta & -\sin \theta\end{pmatrix}$, we have
$$
P\begin{pmatrix} a_{n+1}\cr a_{n+2}\end{pmatrix} = \begin{pmatrix} \cos \theta & -\sin \theta\cr
\sin \theta & \cos \theta\end{pmatrix} P \begin{pmatrix} a_{n}\cr a_{n+1}\end{pmatrix}+ P
\begin{pmatrix} 0\cr \langle \lambda a_{n+1} \rangle\end{pmatrix}
$$
where $\langle x \rangle$ is the fractional part of $x$.
Therefore this gives a rotation map of angle $\theta$ acting on a lattice 
$P \Z^2$ but the image requires 
a bounded perturbation of modulus less than two to fit into lattice points of $P \Z^2$. 
For conjecture \ref{Rotation} we expect that such perturbations do not cumulate and 
the orbits stay bounded, equivalently, all orbits become periodic. 

A nice feature of the map $(a_n,a_{n+1})\mapsto (a_{n+1},a_{n+2})$ is that 
it is clearly bijective on $\Z^2$ by symmetry, while 
under the usual round off scheme, 
the digital information should be more or less lost by the irrational
rotation. This motivates dynamical study of global stability of this algorithm.

The conjecture is trivial when
$\lambda=0,\pm1$. Among non-trivial cases, the second tractable case is 
when $\theta$ is rational and $\lambda$ is quadratic over $\Q$. 
Akiyama, Brunotte, Peth\H o and Steiner \cite{Akiyama-Brunotte-Pethoe-Steiner:07} proved:

\begin{thm}
The conjecture is valid for
$
\lambda=\frac{\pm 1 \pm \sqrt{5}}{2}, \pm \sqrt{2}, \pm \sqrt{3}.
$
\end{thm}

It seems hard 
to prove Conjecture \ref{Rotation}
for other values. 
The case $\lambda=\frac{1-\sqrt{5}}{2}$ was firstly shown 
by Lowenstein, Hatjispyros and Vivaldi \cite{Lowenstein-Hatjispyros-Vivaldi}
with heavy computer assistance. 
A number theoretical proof for
$\frac{1+\sqrt{5}}{2}$ appeared
in \cite{Akiyama-Brunotte-Pethoe-Steiner:06}, whose proof is short but not so easy to
generalize. 
We try to give an accessible 
account using self-inducing piecewise isometry in the
case $\lambda=\omega=\frac{1+\sqrt{5}}{2}$, together with its further dynamical behavior.
The proof in Section \ref{Golden}
is basically in \cite{Akiyama-Brunotte-Pethoe-Steiner:07}.
However this version may elucidate the background idea and is
directly connected to the scaling constant of self-inducing structure
of piecewise isometry acting on a lozenge.

A Pisot number is an algebraic integer $>1$ whose conjugates 
have modulus less than $1$. Throughout the paper, we will see 
the importance of the fact that the scaling constant of self-inducing system
is a Pisot number. Our all discussions heavily depend on this fact. 
Indeed, Pisot scaling constants
often appear in self-inducing structures of several important dynamical systems, 
for e.g., IET and substitutive dynamical systems.
We discuss this point 
in Section \ref{sec:other_self_similar_systems}.
It is pretty surprising that we see 
this phenomenon in cubic piecewise rotations as well. 
We hope this paper gives an easy way to access this interesting area of mathematics. 

We wish to show our gratitude to P.Hubert, 
W.Steiner and F.Vivaldi for helpful comments and relevant literatures
in the development of this manuscript.






\section{Proof of the periodicity for golden mean}
\label{Golden}

Setting $\zeta=\exp(2\pi i/5)$, we have $\omega=-\zeta^{2}-\zeta^{-2}$
and $1/\omega=\zeta+\zeta^{-1}$.
The integer ring of 
$\Q(\zeta)$ coincides with the ring $\Z[\zeta]$
generated by $\zeta$ in $\Z$, $\Z[\zeta]$ is
a free $\Z$-module generated by $1,\zeta,\zeta^2,\zeta^3$.
Hereafter we use a different base as a $\Z$-module:
\begin{lem}
\label{PentagonCord}
$\Z[\zeta]$ is a free $\Z$-module of rank $4$
generated by $1,\omega, \zeta, \omega\zeta$.
\end{lem}

\proof
From 
$\omega=-\zeta^2-\zeta^{-2}$, we have
$$
x_1+x_2 \omega + (y_1+y_2 \omega) \zeta=
(x_1+y_2)+ (y_1+y_2) \zeta + (y_2-x_2) \zeta^2 -x_2 \zeta^3.
$$
On the other hand
$$
a_0+a_1\zeta+a_2\zeta^2+a_3\zeta^3
=
(a_0-a_2+a_3)-a_3 \omega + ((a_1-a_2+a_3)+(a_2-a_3)\omega)\zeta.
$$
\qed
\bigskip

Taking the complex conjugate, the same statement is valid with another 
basis $1,\omega, \zeta^{-1}, \omega\zeta^{-1}$.
Thus 
each element in $\Z[\zeta]$ has a unique expression:
$$
x -\zeta^{-1}y \qquad (x,y\in \Z[\omega]).
$$
Denote by $\langle x \rangle$ the fractional part of $x\in \R$. Then
a small computation gives
\begin{align*}
 0\le a_n +\omega a_{n+1} + a_{n+2} &<1 \\
a_n +\omega a_{n+1} + a_{n+2} & = \langle \omega a_{n+1} \rangle \\
\langle \omega a_n \rangle
- \frac 1{\omega} \langle \omega a_{n+1}\rangle  + \langle \omega a_{n+2} \rangle
& \equiv  0 \pmod{\Z} \\
 x_n - (\zeta +\zeta^{-1}) x_{n+1} + x_{n+2} &\equiv 0 \pmod{\Z}\\
(x_{n+1} - \zeta^{-1} x_{n+2}) & \equiv \zeta^{-1} (x_n - \zeta^{-1} x_{n+1}) \pmod{\zeta^{-1}\Z}
\end{align*}
 and $x_n=\langle \omega a_n \rangle$.
Our problem is therefore embedded into a piecewise isometry $T$ acting on a lozenge 
$[0,1) +(-\zeta^{-1}) [0,1)$:
$$
T(x)=\begin{cases} x/\zeta & \text{Im} (x/\zeta)\ge 0\\
                   (x-1)/\zeta & \text{Im} (x/\zeta)<0
                   \end{cases}.
$$


The action of $T$ is geometrically described in Figure \ref{Piece}.
The lozenge $L=[0,1) +(-\zeta^{-1}) [0,1)$ is rotated by 
the multiplication of $-\zeta^{-1}$ and then the trapezoid $\mathcal{Z}$
which falls outside $L$ is pulled back in by adding $-\zeta^{-1}$.
In total, the isosceles triangle $\Delta$ is rotated clockwise 
by the angle $3\pi/5$ around the origin and the trapezoid $\mathcal{Z}$ is rotated
by the same angle but around the point 
$\frac 12+ i \frac {\sqrt{5(5+2\sqrt{5})}}{10}\simeq 0.5+0.6882 i$ indicated by a black spot, that is
the intersection of two diagonals. 
Our aim is to show that each point $x\in \Z[\zeta] \cap L$ gives 
a periodic $T$-orbit. 

A. Goetz \cite{Goetz:SEPIA} gave a slightly different map. Ours is an 
`inclined' modification of
\cite{Kouptsov-Lowenstein-Vivaldi:02}
and \cite{Akiyama-Brunotte-Pethoe-Steiner:07}.

Clearly the map $T$ is bijective and 
preserves 2-dimensional Lebesgue measure $\mu$. However
the measure dynamical system $(L, \nu, \mathbb{B}, T)$ (with the $\sigma$-algebra 
$\mathbb{B}$ of Lebesgue measurable sets)
is far from ergodic. 
It turned out that orbits of $T$ is periodic for almost all points
but for an exceptional set of Lebesgue measure zero. Our goal is
to prove that the set $\Z[\zeta]$ has no intersection with
this exceptional set. This is not so obvious since 
$\Z[\zeta]$ is dense in $L$ because $\Z[\omega]$ is dense in $\R$.

\begin{figure}[ht]
\begin{center}
\subfigure{
\psfig{figure=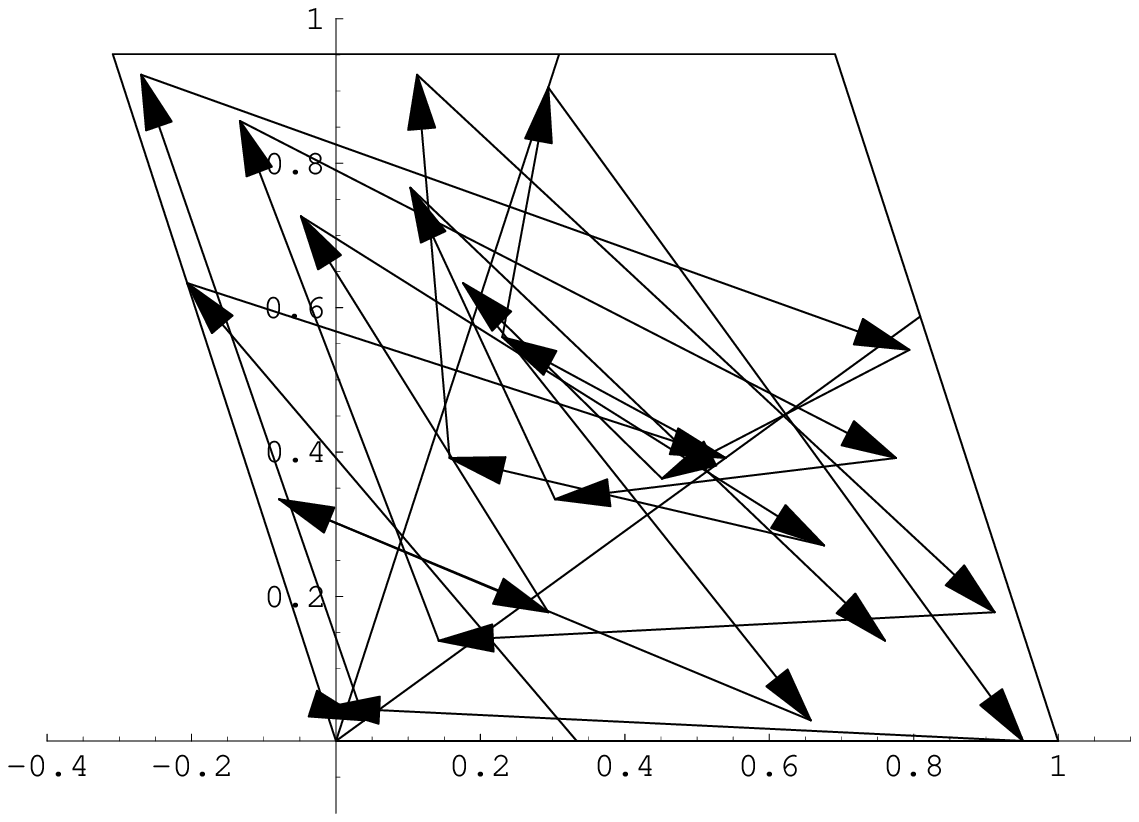, height=4.5cm}
}
\quad
\subfigure{
\psfig{figure=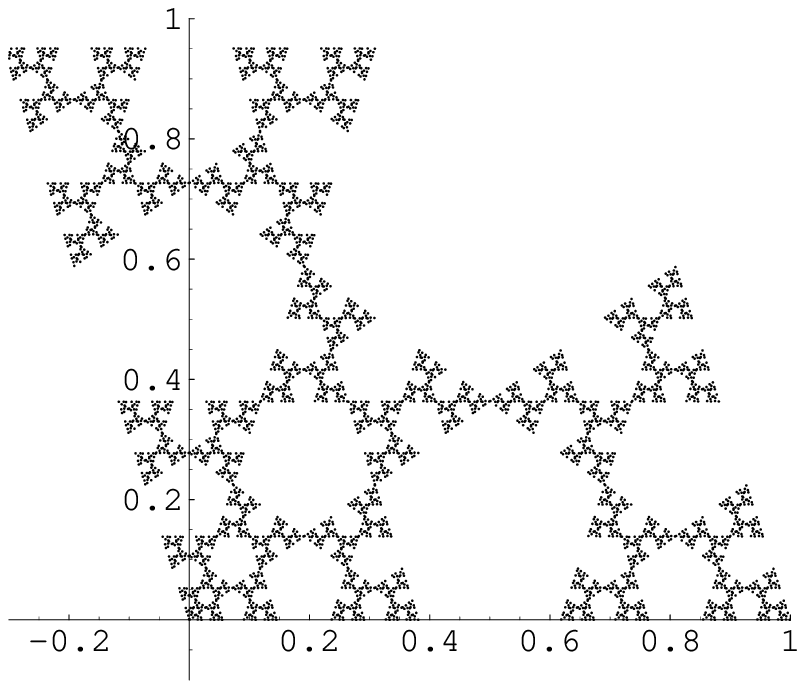, height=4.5cm}
}
\end{center}
\caption{The orbit of $1/3$\label{Orbit2}}
\end{figure}

To illustrate the situation, it is instructive 
to describe an orbit of $1/3$. See Figure \ref{Orbit2}.

Later we will show that the orbit of $1/3$ is aperiodic and forms a
dense subset of the exceptional 
set of aperiodic points. 
Roughly speaking, our task is to show that $\Z[\zeta]\cap L$
has no intersection with the fractal set appeared Figure \ref{Orbit2}. 

\begin{figure}[ht]
\begin{center}
\psfig{figure=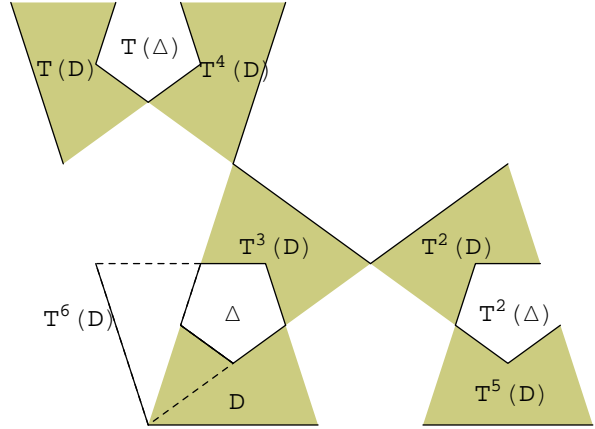}
\end{center}
\caption{Self Inducing structure\label{SelfInducing}}
\end{figure}

The key to the proof is a self-inducing structure with a scaling constant $\omega^2$.
We consider a region $L'=\omega^{-2}L$ and consider the first return map
$$
\hat{T}(x)=T^{m(x)}(x)
$$
for $x\in L'$ where $m(x)$ is the minimum positive integer such that $
T^{m(x)}(x)\in L'$. For any $x\in L'$, the value $m(x)=1,3$ or $6$. 
We can show 
that
\begin{equation}
\label{Self}
\omega^2 \hat{T}(\omega^{-2}x)=T(x)
\end{equation}
for $x\in L$. The proof is geometric, shown in Figure \ref{SelfInducing}. 
The return time $m(x)=3$ in the open pentagon $\Delta'=\omega^{-2}\Delta$ [this is marked $\Delta$ in the figure] and 
$m(x)=6$ in the
shaded pentagonal region $D$ with three closed and two open edges. 
In the remaining isosceles triangle in $\omega^{-2} L$ (whose 
two equal edges are closed and the other open), 
the return time $m(x)$ is $1$.

Note that
the equation is valid for {\bf all} $x\in L'$. 
This makes the later discussion 
very simple. Unfortunately this is not the case for other quadratic values of $\gamma$
and we have to study the behavior of the boundary independently, see
\cite{Akiyama-Brunotte-Pethoe-Steiner:07}. 

\begin{figure}[ht]
\begin{center}
\psfig{figure=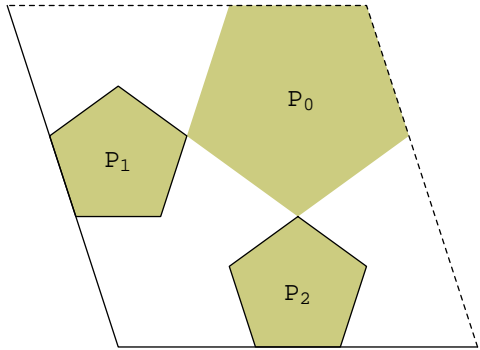}
\caption{Period Pentagons\label{Clopen}}
\end{center}
\end{figure}

Let $U$ be the 1-st hitting map to $L'$ for $x\in L$, i.e., 
$U(x)=T^{m(x)}(x)$ for the minimum non-negative integer $m(x)$ such that 
$T^{m(x)}(x)\in L'$. Note that $U$ is a partial function, i.e., 
$U(x)$ is not defined when there is no positive integer $m$
such that $T^{m}(x)\in L'$.
Since
$$
T(x)=\begin{cases} x/\zeta & x\in \Delta\\
                   (x-1)/\zeta & x\in T(\mathcal{Z})\setminus \Delta
                   \end{cases},
$$
it is easy to make the map $U$ explicit:
$$
U(x)=\begin{cases}
x & x \in L'\\
\left( x-1 \right)/\zeta & x \in T^5(D)\\
\left( x-\zeta \right)/\zeta^2 & x \in T^4(D)\\
\left( x-\frac{\zeta}{\omega} \right)/\zeta^3 & x \in T^3(D)\\
\left(x+\frac{1}{\omega \zeta^2} \right)/\zeta^4 & x \in T^2(D)\\
x+\frac{1}{\omega \zeta} & x \in T(D)\\
\text{Not defined} & x\in P_0\cup P_1 \cup P_2
\end{cases}
$$
where $P_0$ is the largest open pentagon 
and $P_1$ and $P_2=P_1/\zeta$ are two second largest closed 
pentagons in Figure \ref{Clopen}.

Set 
$$
Q=\left\{ 0, 1, \zeta, \frac{\zeta}{\omega}, -\frac{1}{\omega \zeta^2},
-\frac{1}{\omega \zeta} \right\}
=\{d_{0},d_1,d_2,d_3,d_4,d_5\} \subset \Z[\zeta]
$$
to use later. 

We introduce a crucial map $S$ which is the composition of the
1-st hitting map $U$ and expansion by $\omega^2$, i.e. $S(x)=\omega^2 U(x)$.
Denote by $\pi(x)$ the period of $T$-orbits of $x\in L$ and put
$\pi(x)=\infty$ if $x$ is not periodic by $T$. (We
easily see 
$\pi(x)=5$ in $P_0$ and
$\pi(x)=10$ in $P_1 \cup P_2$ unless $x$ is the centroid of the pentagon.)
Then if $\pi(x)$ and $\pi(S(x))$ are defined and finite, 
then we see that $\pi(S(x))<\pi(x)$
which is a consequence of Equation (\ref{Self}).
Therefore if $\pi(x)$ is finite then we have a decreasing sequence
$$
\pi(x)>\pi(S(x))>\pi(S^2(x))>\dots
$$
of positive integers. This shows that there exists a positive integer $k$
such that $S^k(x)$ is not defined. In this case we say that
$S$-orbit of $x$ in finite. We easily see that if $S$-orbit of $x\in L$
is finite, then clearly $\pi(x)$ is finite by Equation (\ref{Self}). 
Thus we have a clear distinction: $x\in L$ is $T$-periodic if and only if
its $S$-orbit is finite.
Assume that $x\in L\cap \Z[\zeta]$ gives an infinite $S$-orbit.
When $U(x)$ is defined, we have
$U(x)=(x-d_{m(x)})/\zeta^{m(x)}$
with $m(x)=\{0,1,2,3,4,5\}$ and $d_i\in Q$ 
for all $x\in L$. 
Thus we have
$$
S^k(x)= \omega^{2k} \frac{x}{\zeta^{\sum_{j=1}^k m_j}}
- \sum_{i=1}^k \omega^{2(k-i+1)} \frac{d_{m_i}}{\zeta^{\sum_{j=i}^{k} m_j}}.
$$
By the assumption $S^k(x)$ is defined for $k=1,2,\dots$
and stays in $L$. 
Consider the conjugate map $\phi$ which sends $\zeta\rightarrow \zeta^2$.
As $\phi(\omega)=-1/\omega$, we have
$$
\phi(S^k(x))= \frac {\phi(x)}{\omega^{2k} \zeta^{2\sum_{j=1}^k m_j}}
- \sum_{i=1}^k \frac{d'_{m_i}}{\omega^{2(k-i+1)} \zeta^{2\sum_{j=i}^{k} m_j}}
$$
with $d'_i=\phi(d_i)\in \phi(Q)$. 
Put $A=\max\{|d'_i|\ :\ d_i\in Q\}$. Then we have
$$
|\phi(S^k(x))| \le |\phi(x)| + \frac A{\omega^{2}-1}
$$
Thus we have $S^k(x), \phi(S^k(x))$ and their complex conjugates are bounded
by a constant which does not depend on $k$. This implies that the sequence
$(S^k(x))_{k}$ must be eventually periodic.

Summing up, for a point $x$ in $\Z[\zeta]$, its $S$-orbit is finite or
eventually periodic. When it is finite then its $T$-orbit is periodic and
when its $S$-orbit is eventually periodic then $T$-orbit is aperiodic.
Thus we have an algorithm for $x\in \Z[\beta]\cap L$ to tell
whether its $T$-orbit is periodic or not.
Since 
$$
|\phi(S^k(x))| \le \frac{|\phi(x)|}{\omega^{2k}} + \frac A{1-\omega^{-2}},
$$
for any positive $\varepsilon$, 
the right hand side is bounded by
$$
\varepsilon+\frac A{\omega^{2}-1}
$$
for a sufficiently large $k$. This means that under the assumption that
there is an infinite $S$-orbit, the set
$$
\left\{x\in \Z[\zeta]\cap L\ :\ |\phi(x)|\le \varepsilon+\frac A{\omega^{2}-1}\right\}
$$
contains $x$ with $\pi(x)=\infty$. Since this set is finite, 
it is equal to 
$$
B=\left\{x\in \Z[\zeta]\cap L\ :\ |\phi(x)|\le \frac A{\omega^{2}-1}\right\}
$$
for a sufficiently small $\varepsilon$. 
Since there are only finitely many candidates in $B$,
we obtain an algorithm to check
whether an element $x\in \Z[\zeta]\cap L$ with $\pi(x)=\infty$ exists.
In fact, all elements in $B$ gives a finite $S$-expansion,
we are done. 

The same algorithm applies to $\frac 1M \Z[\zeta]$ 
with a fixed positive integer $M$. In this way, 
we can also show that points in $\frac 12 \Z[\zeta]$ are periodic.
We can find aperiodic orbits in $\frac 13 \Z[\zeta]$. 
For example, one can see that $1/3$ has 
an aperiodic $T$-orbit because its $S$-orbit:
$$
\frac13,\frac{w^2}3, -\frac{\zeta^{-1}}3, -\frac{\omega^2 \zeta^{-1}}3
-\frac{2\zeta^{-1}}3, -\frac{\omega^{-2}\zeta^{-1}}3, -\frac{\zeta^{-1}}3, \dots
$$
satisfies $S^2(1/3)=S^6(1/3)$.

It is crucial in the above proof that 
the scaling constant of the self-inducing structure is a Pisot number. 
Scaling constants of piecewise isometries often become Pisot numbers, moreover 
algebraic units. We discuss these phenomena in Section
\ref{sec:other_self_similar_systems}.

\section{Coding of aperiodic $T$-orbits}

\begin{figure}[ht]
\begin{center}
\subfigure[All digits\label{Full}]{
\psfig{figure=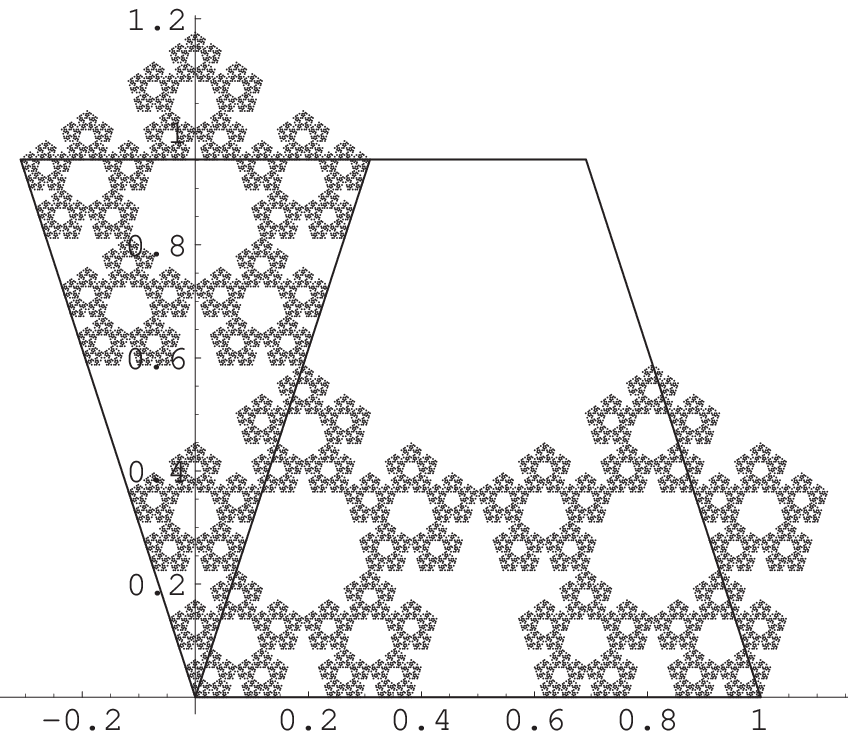, height=5cm}
}
\quad
\subfigure[$d_0,d_2,d_3,d_5$\label{Essential}]{
\psfig{figure=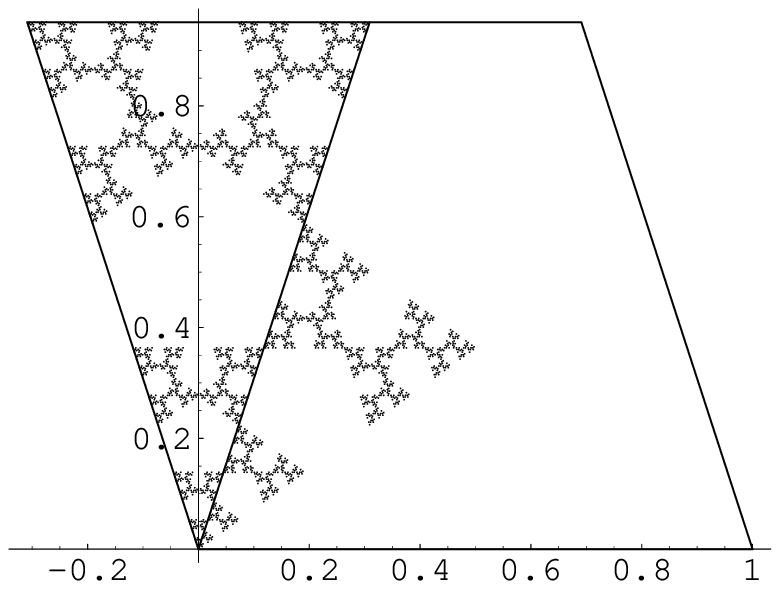, height=4cm}
}
\caption{Attractors containing $\mathbf{A}$}
\end{center}
\end{figure}

Denote by $\mathbf{A}$ the set of all $T$-aperiodic points in $L$. 
By the proof of the previous section, we have
$$
\mathbf{A}=\{ x\in L\ |\ S^k(x) \text{ is defined for all } k=1,2,\dots \}. 
$$
We also have $S(\mathbf{A})\subset \mathbf{A}$. 
This means that for
$x_1\in \mathbf{A}$, there is a $m_i\in \{0,1,2,3,4,5\}$ and $x_i\in \mathbf{A}$ 
such that $\omega^2\zeta^{-m_i}(x_i-d_{m_i})=x_{i+1}\in \mathbf{A}$ 
for $i=1,2,\dots$. We therefore
have an expansion
\begin{equation}
\label{beta}
x_1=d_{m_1}
+\frac {\zeta^{m_1}}{\omega^2}\left(
d_{m_2}+\frac {\zeta^{m_2}}{\omega^2}\left(
d_{m_3}+\frac {\zeta^{m_3}}{\omega^2}\left(
d_{m_4}+\frac {\zeta^{m_4}}{\omega^2}
\dots \right.\right.\right.
\end{equation}
Conversely a sequence $\{m_i\}_{i=1,2,\dots}$ defines a single point of $Y$. 
Therefore $\mathbf{A}$ must be a subset of the attractor $Y$ 
of the iterated function system (IFS):
$$
Y=\bigcup_{i=0}^5 \left(\frac {\zeta^i}{\omega^2} Y+d_i\right),
$$
an approximation of which is depicted in Figure \ref{Full}.
At this point we can assert that $2$-dimensional
Lebesgue measure of aperiodic points in $L$ must be 
zero, because $\omega^4\simeq 6.854\dots>6$. 

We notice that the digits in $Q$ are not arbitrarily chosen because the
image of $S$ must be in $T(\mathcal{Z})$. Thus the digits $d_1$ and $d_4$
appears only at 
the beginning in the expression of Equation (\ref{beta}). Therefore
it is more suitable to study $\mathbf{A}\cap T(\mathcal{Z})$. 
The attractor
\begin{equation}
\label{attractor}
Y'=\left(\frac 1{\omega^2} Y'+d_0\right) \cup 
\left(\frac{\zeta^2}{\omega^2} Y'+d_2 \right)\cup
\left(\frac{\zeta^3}{\omega^2} Y'+d_3\right) \cup
\left(\frac{\zeta^5}{\omega^2} Y'+d_5\right)
\end{equation}
is depicted in Figure \ref{Essential}. 

\begin{figure}[ht]
\begin{center}
\subfigure{
\psfig{figure=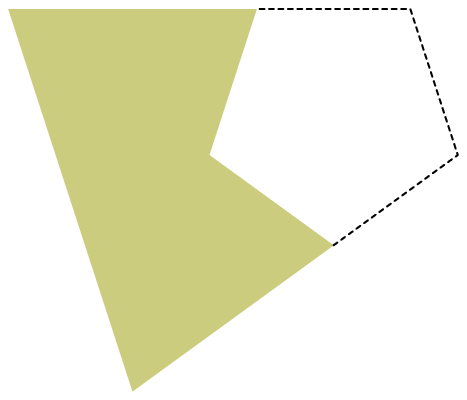}
}
\quad
\subfigure{
\psfig{figure=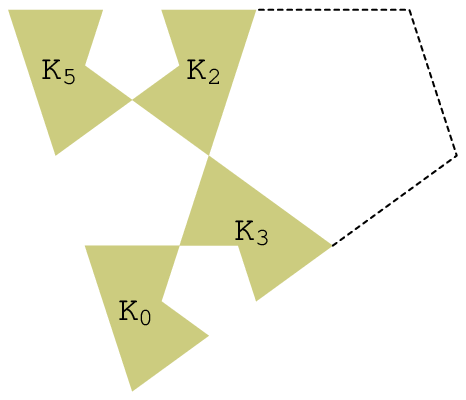}
}
\end{center}
\caption{Open set condition\label{OSC}}
\end{figure}

This iterated function system satisfies
OSC by a pentagonal shape $K$ with whose vertices are
$$
0,-\zeta^{-1}, \zeta, 
-\zeta \omega^{-1}-\zeta^{-1},
-\zeta^2 \omega^{-1}
$$
as in Figure \ref{OSC}. 
We confirm that 
the pieces $K_m=\frac {\zeta^m}{\omega^2} K +d_m$ do not 
overlap. 

\begin{figure}[ht]
\begin{center}
\subfigure{
\psfig{figure=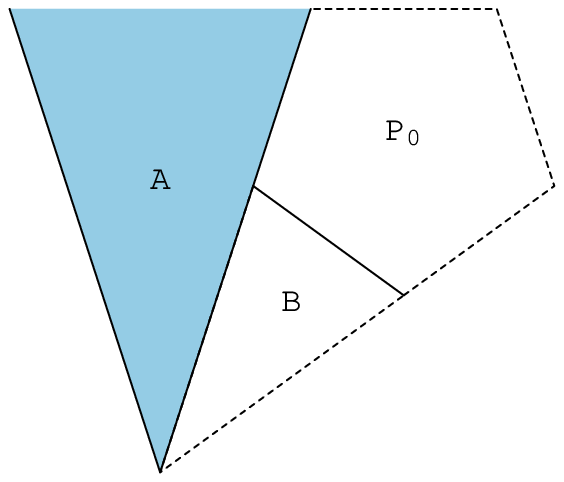}
}
\quad
\subfigure{
\psfig{figure=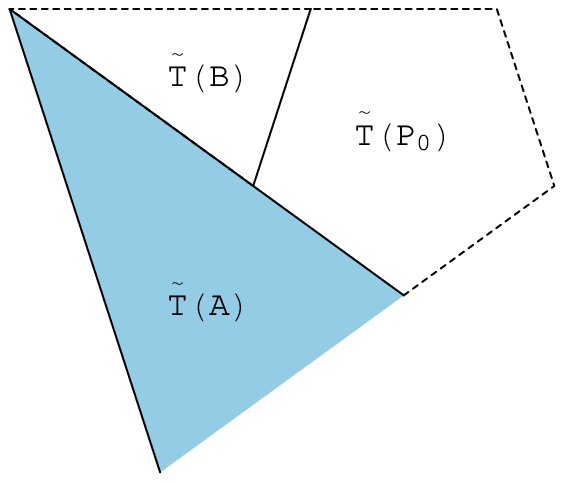}
}
\end{center}
\caption{Induced Rotation $\widetilde{T}$ 
on $T(\mathcal{Z})$\label{InducedRotation}}
\end{figure}

We consider the induced system of $(L,\mathbb{B},\nu,T)$ to
$T(\mathcal{Z})$.
Denote by $\widetilde{T}$ the first return map 
on $T(\mathcal{Z})$. Then the induced system $(T(\mathcal{Z}),
\widetilde{T})$ is the domain exchange of two isosceles triangle $A$ and $B$ depicted in
Figure \ref{InducedRotation}. The triangle $A$ has two closed edges of equal length 
and one open edge,
while $B$ has one closed edge and two open edges of the same length. 
The open regular pentagon 
$P_0$ and the triangle $B$ move together by $\widetilde{T}$ 
and can be merged into a single shape. 

We see
\begin{equation}
\label{Induction0}
\widetilde{T}(x)= \begin{cases} T^2(x) & x\in \Delta\\
                      T(x) & x\in T(\mathcal{Z})\setminus \Delta.
\end{cases}
\end{equation}

\begin{figure}[ht]
\begin{center}
\psfig{figure=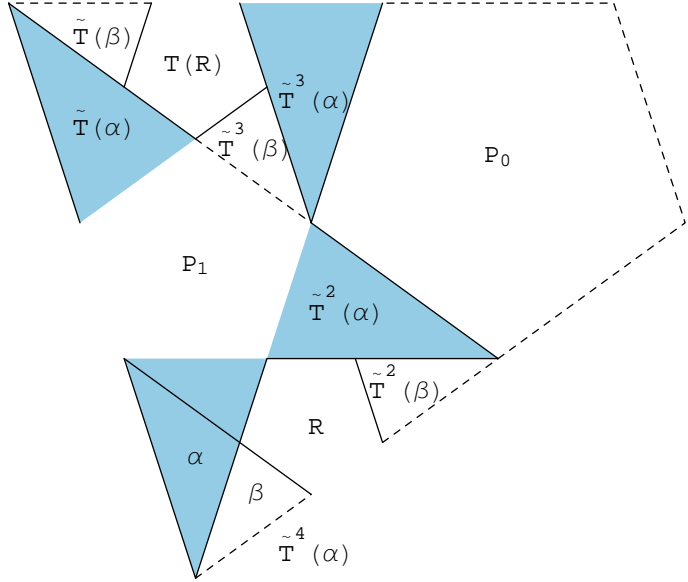}
\caption{Self Inducing Structure of $(T(\mathcal{Z}),\widetilde{T})$\label{SelfInducing1}}
\end{center}
\end{figure}

Again we find self-inducing structure with the scaling 
constant $\omega^2$: 
\begin{equation}
\label{Self2}
\omega^2 \widetilde{T}(\omega^{-2} x)= \widetilde{T}(x)
\end{equation}
for all $x\in T(\mathcal{Z})$. This can be seen 
in Figure \ref{SelfInducing1} with 
$\alpha=\omega^{-2}A$, $\beta=\omega^{-2} B$ 
and $R=\omega^{-2}P_0$.
This induced 
dynamics $(T(\mathcal{Z}), \widetilde{T})$ is essential in describing the set 
$\mathbf{A}$.

Readers may notice that we can find a 
self-inducing structure by smaller scaling constant $\omega$ in 
Figure \ref{Essential} by taking two connected pieces. 
However this choice of inducing region 
is not suitable because the self-inducing relation (with flipping)
is measure theoretically valid, but has different behavior on the boundary. 

Let us introduce two codings. 
First is the coding of $T$-orbits of a point $x$ 
in $L$ in two symbols $\{0,1\}$: 
$\mathbf{d}(x)=(\psi(T^n(x))_n\in \{0,1\}^{\N}$ where 
$$\psi(x)=\begin{cases} 0 & x\in \Delta\\
                              1 & x\in \mathcal{Z}
\end{cases}.$$
For e.g.,
the $\mathbf{d}(1/3)=10110101011010101101101101\dots$ The second coding is defined by
$\widetilde{\mathbf{d}}(x)= (\widetilde{\psi}(\widetilde{T}^n(x)))_n \in \{a,b\}^{\N}$
for $x\in T(\mathcal{Z})$ where 
$$\widetilde{\psi}(x)=\begin{cases} a & x\in \Delta\\
                           b & x\in T(\mathcal{Z})\setminus \Delta.
\end{cases}$$
For a point $x$ in $T(\mathcal{Z})$ we have two codings by $\{a,b\}$ and
by $\{0,1\}$. From Equation (\ref{Induction0}), the two codings are equivalent through the 
substitution $a\rightarrow 01$, $b\rightarrow 1$. For a given coding
of $T$-orbit by $\{0,1\}$, there is a unique way to retrieve the coding
of $\widetilde{T}$-orbit by $\{a,b\}$, because the symbol $0$ must be followed by $1$.
For e.g., $T(1/3)=-2\zeta^{-1}/3\in T(\mathcal{Z})$ is coded in two ways as:
$$
\widetilde{\mathbf{d}}(-2\zeta^{-1}/3)=
a\: b\: a\: a\: a\: b\: a\: a\: a\: b\: a\: b\: a\: b\: a\: b\dots
$$
and
$$
\mathbf{d}(-2\zeta^{-1}/3)=
01\:1\:01\:01\:01\:1\:01\:01\:01\:1\:01\:1\:01\:1\:01\:1\dots
$$
Hereafter we discuss the coding $\widetilde{\mathbf{d}}$.
Observing the trajectory of  
the region $\omega^{-2}(\Delta)$ and
$\omega^{-2}(T(\mathcal{Z})\setminus \Delta)$ 
 by the first return map by the iteration of $\widetilde{T}$ to 
the region $\omega^{-2}T(\mathcal{Z})$, it is natural to introduce 
a substitution $\sigma_0$:
$$
a\rightarrow aaba,\quad b\rightarrow baba.
$$
on $\{a,b\}^*$ and we have
$$
\widetilde{\mathbf{d}}(\omega^{-2}x)=\sigma_0(\widetilde{\mathbf{d}}(x))
$$
for $x\in T(\mathcal{Z})$. More generally, following the analogy of the previous section, the first hitting map 
to the region $\omega^{-2}(T(\mathcal{Z}))$ provide us an expansion of a point $x\in 
T(\mathcal{Z})$ exactly
in the same form as (\ref{beta}) with restricted digits $\{d_0,d_2,d_3,d_5\}$.
One can confirm that
\begin{equation}
\label{Delay}
\widetilde{\mathbf{d}}\left( \frac{\zeta^m}{\omega^2} x + d_m\right)=\begin{cases}
\sigma_0(\widetilde{\mathbf{d}}(x)) & m=0\\
a \oplus \sigma_0(\widetilde{\mathbf{d}}(\widetilde{T}(x))) & m=2\\
ba \oplus \sigma_0(\widetilde{\mathbf{d}}(\widetilde{T}^2(x)))& m=3\\
aba \oplus \sigma_0(\widetilde{\mathbf{d}}(\widetilde{T}^3(x)))& m=5
\end{cases}
\end{equation}
where $\oplus$ is the concatenation of letters. Defining
conjugate substitutions by
$\sigma_1=a\sigma_0 a^{-1}$, $\sigma_2=ba\sigma_1 a^{-1}b^{-1}$ and
$\sigma_3=aba \sigma_0 a^{-1}b^{-1}a^{-1}$, i.e.,
\begin{eqnarray*}
\label{Sadic}
\sigma_0(a)=aaba, &&\quad \sigma_0(b)=baba\\
\sigma_1(a)=aaab, &&\quad \sigma_1(b)=abab\\
\sigma_2(a)=baaa, &&\quad \sigma_2(b)=baba\\
\sigma_3(a)=abaa, &&\quad \sigma_3(b)=abab
\end{eqnarray*}
one may rewrite
$$
\widetilde{\mathbf{d}}\left( \frac{\zeta^m}{\omega^2} x + d_m\right)=\begin{cases}
\sigma_0(\widetilde{\mathbf{d}}(x)) & m=0\\
\sigma_1(\widetilde{\mathbf{d}}(\widetilde{T}(x))) & m=2\\
\sigma_2(\widetilde{\mathbf{d}}(\widetilde{T}^2(x)))& m=3\\
\sigma_3(\widetilde{\mathbf{d}}(\widetilde{T}^3(x)))& m=5.
\end{cases}
$$
We say that 
an infinite word $y$ in $\{a,b\}^{\N}$ is an $S$-adic limit of 
$\sigma_{i}\ (i=0,1,2,3)$ 
if there exist $y_i\in \{a,b\}^{\N}$ for $i=1,2,\dots$ such that
$$
y=\lim_{\ell\rightarrow \infty}
\sigma_{m_1}\circ \sigma_{m_2} \circ \sigma_{m_3}\circ \cdots \circ 
\sigma_{m_{\ell}}(y_{\ell}).
$$
with $m_i\in \{0,1,2,3\}$.
Since each element $x\in T(\mathcal{Z})\cap \mathbf{A}$ has an infinite expansion
(\ref{beta}) with digits $\{d_0,d_2,d_3,d_5\}$, we find 
$x_i\in \omega^{-2}T(\mathcal{Z})$
such that
$$
\widetilde{\mathbf{d}}(x)=\lim_{\ell\rightarrow \infty}
\sigma_{m_1}\circ \sigma_{m_2} \circ \sigma_{m_3}\circ \cdots \circ
\sigma_{m_{\ell}}(\widetilde{\mathbf{d}}(x_{\ell})).
$$
This shows that $\widetilde{\mathbf{d}}(x)$ is an $S$-adic limit
of $\sigma_{i}\ (i=0,1,2,3)$. 

Note that from the definition (\ref{Sadic}) of $\sigma_i$, 
for a given $S$-adic limit $y$ 
there is an algorithm to retrieve uniquely 
the sequence $(\sigma_{m_i})_i$. 
Checking first four letters of $y$,
we know the first letter of $y_1$ and to determine $m_1$ we need first
6 letters. We can iterate this process easily.

Summing up,
we embedded the set $\mathbf{A} \cap T(\mathcal{Z})$ into the attractor $Y'$ 
of an IFS (\ref{attractor}) and succeeded in characterizing the coding of $\widetilde{T}$-orbits
of points in this attractor as a set of 
$S$-adic limits on $\{\sigma_0,\sigma_1,\sigma_2,\sigma_3\}$.
However recalling that points in closed pentagons $P_1$ and $P_2$ are 
$T$-periodic and $Y'$ is a non-empty compact set, 
we see from Figure \ref{Essential} 
that $\mathbf{A}$ is a proper subset of $Y'$. 

We wish to characterize the set of aperiodic points
in $Y'$ and its coding through $\widetilde{\mathbf{d}}$. Recalling
the discussion in the previous section, if $x\in T(\mathcal{Z})$ has periodic 
$T$-orbits if and only if there exists a positive integer $k$ such that
$S^k(x)\in P_0\cup P_1 \cup P_2$. The equivalent statement in the induced system 
$(T(\mathcal{Z}), \widetilde{T})$ is that $x\in T(\mathcal{Z})$ is
$\widetilde{T}$-periodic if and only if there exists a positive integer
$k$ such that $S^k(x)\in P_0\cup P_1$. Note that we have:
$$
\widetilde{T}(x)=\begin{cases} \zeta^{-1} (x-p)+p & x\in P_0\\
                               \zeta^{-2} (x-q)+q & x\in P_1\\
\end{cases}
$$
where $p=\frac 12+ i \frac {\sqrt{5(5+2\sqrt{5})}}{10}$ 
(resp. $q=i\sqrt{\frac {5+\sqrt{5}}{10}}$) 
is the center of $P_0$ (resp. $P_1$) and consequently
$\widetilde{T}^5(x)=x$ holds for $x\in P_0\cup P_1$. If $x\in T(\mathcal{Z})$ and $x$ is 
$\widetilde{T}$-periodic, then there exist $x_i\in T(\mathcal{Z})$ such that 
$x_{\ell}\in P_0\cup P_1$ and
$$
x=d_{m_1}
+\frac {\zeta^{m_1}}{\omega^2}\left(
d_{m_2}+\frac {\zeta^{m_2}}{\omega^2}\left(
d_{m_3}+\frac {\zeta^{m_3}}{\omega^2}
\dots 
\left(d_{m_{\ell}}+\frac {\zeta^{m_{\ell}}}{\omega^{2\ell}}x_{\ell}
\right) \dots 
\right)\right),
$$
with $m_i\in \{0,2,3,5\}$. Thus the set of $\widetilde{T}$-periodic points in $T(\mathcal{Z})$ 
consists of all the pentagons of the form
\begin{equation}
\label{Pent}
d_{m_1}
+\frac {\zeta^{m_1}}{\omega^2}\left(
d_{m_2}+\frac {\zeta^{m_2}}{\omega^2}\left(
d_{m_3}+\frac {\zeta^{m_3}}{\omega^2}
\dots \left(d_{m_{\ell}}+\frac {\zeta^{m_{\ell}}}{\omega^{2\ell}}P_j
\right) \dots
\right)\right)
\end{equation}
with $j=0,1$ and $m_j\in \{0,2,3,5\}$.
From the self-inducing structure (\ref{Self2}), it is easy to see that
if two points $x,x'$ are in the same pentagon of above shape and none of them
is the center, then they have 
exactly the same periods. Moreover two $\widetilde{T}$-orbits 
keeps constant distance, i.e.,
$\widetilde{T}^n(x)-\widetilde{T}^n(x')=\zeta^s (x-x')$ for some integer $s$.
The period is completely determined in
\cite{Akiyama-Brunotte-Pethoe-Steiner:07}. 
We have all the periodic orbits in $T(\mathcal{Z})$ and therefore have 
a geometric description of aperiodic points: 
$$
\mathbf{A} \cap T(\mathcal{Z}) = 
T(\mathcal{Z}) \setminus \{\text{All pentagons of the form (\ref{Pent})} \}.
$$

\begin{figure}[ht]
\begin{center}
\subfigure{
\psfig{figure=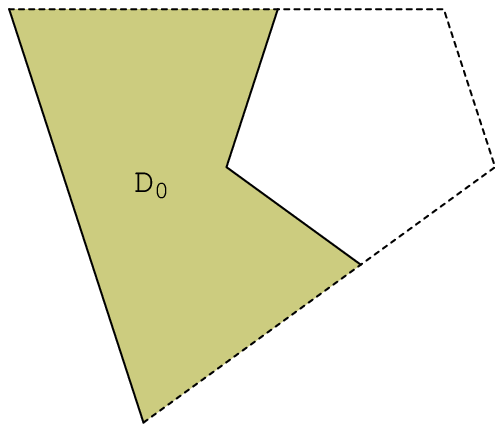}
}
\quad
\subfigure{
\psfig{figure=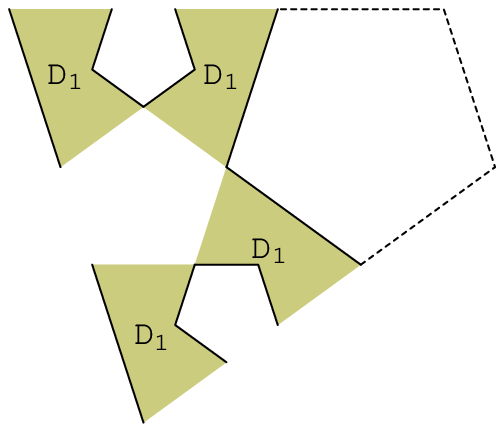}
}
\end{center}
\caption{Pentagon Removal Algorithm\label{Delete}}
\end{figure}

\begin{figure}[ht]
\begin{center}
\psfig{figure=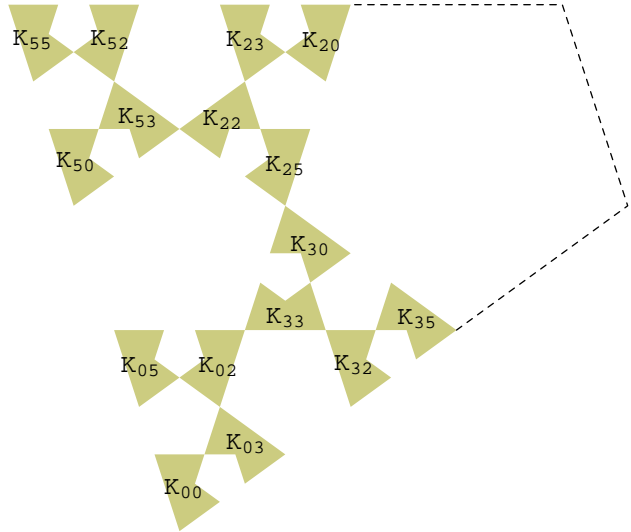}
\end{center}
\caption{Subdivision procedure\label{OSC2}}
\end{figure}

Subtraction of these pentagons from $T(\mathcal{Z})$ is described 
by an algorithm. The initial set is 
$D_0=T(\mathcal{Z})\setminus P_0$ with two open and three closed edges as
 in the left Figure \ref{Delete}. The interior ${\rm Inn}(D_0)$ gives another feasible
 open set to assure the open set condition of the IFS of (\ref{attractor}).
Inductively we define the decreasing sequence of sets
$$
D_{i+1}= \bigcup_{m\in \{0,2,3,5\}} \left(\frac {\zeta^m}{\omega^2} D_i + d_m
\right)
$$
for $i=0,2,\dots$. Then $D_i$ consists of $4^i$ pieces congruent
to $\omega^{-2i} D_0$ without overlapping. Note that since
$$
D_1=D_0 \setminus (P_1 \cup \omega^{-2}P_0 \cup \frac {\omega^{-2}P_0-1}{\zeta}),
$$
$D_1$ is obtained by subtracting from $D_0$ one closed and two open regular pentagons
as in Figure \ref{Delete}. 
To generate $D_{i+1}$, 
each $4^i$ pieces in $D_i$ are subdivided into $4$ sub-pieces by subtracting
three small regular pentagons. Clearly all regular pentagons
of the shape (\ref{Pent}) are subtracted by this iteration and we obtain
$$
\mathbf{A} \cap T(\mathcal{Z}) = \bigcap_{i=0}^{\infty} D_i.
$$
This observation allows us to symbolically characterize aperiodic points
in $Y'$. 
First, every point $x$ of $Y'$ has an address
$d_{m_1} d_{m_2} \dots \in \{d_0,d_2,d_3,d_5\}^{\N}$ by the expansion (\ref{beta}). 
The address is unique but for countable exceptions. The exceptional
points forms the set of cut points of $Y'$ having the eventually periodic
expansion:
\begin{eqnarray*}
d_0 d_2 (d_0)^{\infty} &\simeq& d_3 d_3 (d_5)^{\infty} \\
d_3 (d_0)^{\infty} &\simeq& d_2 (d_5)^{\infty} \\
d_2 d_2 (d_0)^{\infty} &\simeq& d_5 d_3 (d_5)^{\infty}
\end{eqnarray*}
in the suffix of its address, which is understood by Figure \ref{OSC2} where
$K_{mn}=\frac {\zeta^m}{\omega^2} (\frac {\zeta^n}{\omega^2}K+d_n) +d_m$.

Note that if a point $x$ in $T(\mathcal{Z})$ 
is periodic, then there exists a non-negative integer
$k$ such that $S^{k}(x)\in P_0 \cup P_1$. Moreover, if $x\in Y'\cap T(\mathcal{Z})$,
then there exists a non-negative integer
$k$ such that $S^{k}(x)\in \partial(P_1)$, because it can not be 
an inner point of $P_0$ or $P_1$. In other words, such $x$ must be
located in the open edge of one of $4^k$ pieces of $D_k$. 
From Figure \ref{Delete}, one can construct the following Figure \ref{OpenEdge}
which recognize points  of two open edges in $\partial(D_0)$.
For construction, we introduce a new symbol set $\{R,L\}$ (right and left) 
to distinguish which open edge of $D_k$ is into focus. 


To read the graph and obtain the previous sequences, ignore $\{R,L\}$ and substitute $\{0,2,3,5\}$ with $\{d_0,d_2,d_3,d_5\}$. 
A point $x\in Y'$ is periodic (or in the open edge of $D_0$) 
if and only if a suffix of
the address $d_{m_1} d_{m_2} \dots \in \{d_0,d_2,d_3,d_5\}^{\N}$ is 
in Figure \ref{OpenEdge}. 
Note that the points with double addresses
are on the open edge of some $D_i$ and consequently 
their suffixes are read in Figure \ref{OpenEdge}. 
Figure \ref{OSC2}
helps this construction. For e.g., the right open edge of $K_5$ consists of 
the left open edge of $K_{53}$ and the right open edge of $K_{50}$, therefore
we draw outgoing edges from $5R$ to $3L$ and $0R$. 

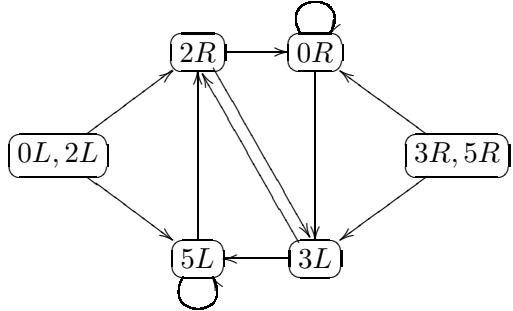
\begin{figure}
\[
\xymatrix{
  & *+[F-:<4pt>]{2R} \ar@<0.5ex>[ddr] \ar[r] &  *+[F-:<4pt>]{0R} \ar@(ul,ur) \ar[dd] & \\
*+[F-:<4pt>]{0L,2L}\ar[dr] \ar[ur] &   &  & *+[F-:<4pt>]{3R,5R} \ar[dl] \ar[ul]\\
  & *+[F-:<4pt>]{5L} \ar@(dl,dr) \ar[uu] & *+[F-:<4pt>]{3L} \ar[l] \ar@<0.5ex>[uul] & \\
}
\]
\caption{$\widetilde{T}$-periodic expansions\label{OpenEdge}}
\end{figure}

\begin{figure}
\[
\xymatrix{
  \ar[rr]&                & *++[o][F]{} \ar@/^/[ddr]^{3,5} \ar@/_2pc/[ddll]_{0,2} \ar@(ul,ur)^{0,2,3,5}  & \\
  & *++[o][F=]{} \ar@<0.5ex>[ddr]^3 \ar[r]^0 &  *++[o][F=]{} \ar@(ul,ur)^0 \ar[dd]^3 & \\
*++[o][F]{} \ar[dr]^5 \ar[ur]^2 &   &  & *++[o][F]{} \ar[dl]^3 \ar[ul]^0\\
  & *++[o][F=]{} \ar@(dl,dr)_5 \ar[uu]^2 & *++[o][F=]{} \ar[l]^5 \ar@<0.5ex>[uul]^2 & 
  \\
}
\]
\caption{B\"uchi automaton for periodic points in $Y'$ \label{OpenEdge2}}
\end{figure}
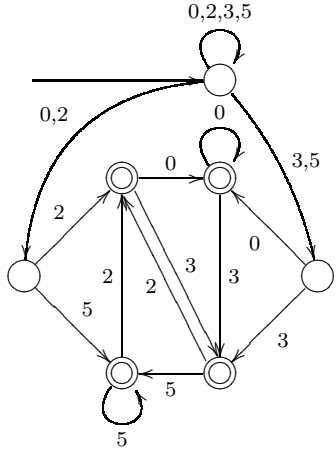

As a result, the set of addresses of the
points in $\mathbf{A} \cap T(\mathcal{Z})$ are recognized by a B\"uchi automaton
which is the complement of the B\"uchi automaton of Figure \ref{OpenEdge2}. Here
the double bordered states in Figure \ref{OpenEdge2} are final states.
Each infinite word produced by the edge labels $\{d_0,d_2,d_3,d_5\}$ 
on this directed graph is accepted,
because it visits infinitely many times the final states.
We do not give here the exact shape of its complement. It is known that 
complementation of a B\"uchi automaton is much harder than the one 
of a finite automaton, because the subset construction 
does not work (c.f. \cite{Thomas:90, Perrin-Pin:04}). 

Now consider the topology of $\{a,b\}^{\N}$ induced from the
metric defined by $2^{-\max_{x_i\neq y_i}i}$ for
$x=x_1x_2\dots, y=y_1y_2\dots \in \{a,b\}^{\N}$.
Take 
a fixed point $w=(w_i)_{i=0,1,2,\dots}\in \{a,b\}^{\N}$ 
with $\sigma_0(w)=w$. This is computed for e.g., by $\lim_n \sigma_0^n(a)$.
The shift map $V$ is a 
continuous map from $\{a,b\}^{\N}$ to itself
defined by $V((w_i))=(w_{i+1})$. Letting  $X_{\sigma_0}$ be the closure of the set 
$\{ V^n (w) \ |\ n=0,1,\dots \}$, we can define the substitutive dynamical system $(X_{\sigma_0}, V)$ associated with $\sigma_0$. 
Since $\sigma_0$ is primitive
the set $X_{\sigma_0}$ does not depend on the choice of the fixed point
and $(X_{\sigma_0}, V)$ is minimal and uniquely ergodic (see \cite{Fogg:02}). 
Let $\tau$ be the invariant measure of $(X_{\sigma_0}, V)$.
On the other hand, for the attractor $Y'$ 
there is the self-similar measure $\nu$, i.e., a unique probability measure
(c.f. Hutchinson \cite{Hutchinson:81})
satisfying
$$
\nu(X)=\frac 14\sum_{m\in \{0,2,3,5\}} 
\nu \left( \frac{\omega^2}{\zeta^m} (X -d_m) \right)
$$
for $\nu$-measurable sets $\mathbb{B}_{Y'}$ in $Y'$. 

\begin{thm}
\label{Main}
The restriction of $\widetilde{T}$ to $Y'$ 
is measure preserving and 
$(Y',\mathbb{B}_{Y'}, \nu, \widetilde{T})$ 
is isomorphic to the $2$-adic odometer $(\Z_2, x\mapsto x+1)$
as measure dynamical systems:
\begin{equation}
\label{Additive}
\begin{CD}
\Z_2  @>+1>> \Z_2\\
@V{\phi}VV              @V{\phi}VV\\
Y' @>\widetilde{T}>> Y'
\end{CD}
\end{equation}
where $\phi:\Z_2 \rightarrow Y'$ is almost one to one and 
measure preserving, which will be made explicit in the proof. 
Moreover the map 
$$
\rho: x\mapsto \frac {x-(x \bmod{4})}{4} 
$$
from $\Z_2$ to itself gives a commutative diagram:
\begin{equation}
\label{Multiplicative}
\begin{CD}
\Z_2  @>\rho>> \Z_2\\
@V{\phi}VV              @V{\phi}VV\\
Y' @>S>> Y'.
\end{CD}
\end{equation}
\end{thm}

The above theorem may be read that 
$(Y',\mathbb{B}_{Y'}, \nu, \widetilde{T})$ 
gives a one-sided variant of numeration system in the sense of 
Kamae \cite{Kamae:05}.

\proof
First we confirm that $\widetilde{T}$ is measure preserving.
Denote by 
$[d_{m_1},d_{m_2},\dots,d_{m_{\ell}}]$ the cylinder set:
\begin{equation}
\label{Cylinder}
d_{m_1}
+\frac {\zeta^{m_1}}{\omega^2}\left(
d_{m_2}+\frac {\zeta^{m_2}}{\omega^2}\left(
d_{m_3}+\frac {\zeta^{m_3}}{\omega^2}
\dots \left(d_{m_{\ell}}+\frac {\zeta^{m_{\ell}}}{\omega^{2\ell}}Y'
\right) \dots
\right)\right)
\end{equation}
By the OSC, we have $\nu([d_{m_1},d_{m_2},\dots,d_{m_{\ell}}])=4^{-\ell}$.
From Figure \ref{SelfInducing1}, we see that 
$\widetilde{T}^{-1}([d_3])=[d_5]$,
$\widetilde{T}^{-1}([d_2])=[d_3]$,
$\widetilde{T}^{-1}([d_0])=[d_2]$ but
$\widetilde{T}^{-1}([d_5])$ intersects both $A$ and $B$.
Hence if $m_1=0,2,3$, then
$\nu(\widetilde{T}^{-1}([d_{m_1},d_{m_2},\dots,d_{m_{\ell}}]))
=4^{-\ell}$. By using the self-inducing structure in Figure \ref{SelfInducing1}, 
we also have
$\widetilde{T}^{-1}([d_5d_3])=[d_0d_5]$,
$\widetilde{T}^{-1}([d_5d_2])=[d_0d_3]$ and
$\widetilde{T}^{-1}([d_5d_0])=[d_0d_2]$.
Thus if $m_2=0,2,3$, then
$\nu(\widetilde{T}^{-1}([d_{5},d_{m_2},\dots,d_{m_{\ell}}]))
=4^{-\ell}$. Repeating this, we can show that
$$\nu(\widetilde{T}^{-1}([d_{m_1},d_{m_2},\dots,d_{m_{\ell}}]))
=4^{-\ell}$$ holds for all $m_i\in \{0,2,3,5\}$ but a single exception
$m_1=m_2=\dots=m_{\ell}=5$. Since $\ell$ is arbitrary chosen, 
a simple approximation argument shows that
$\widetilde{T}$ is measure preserving and 
$(Y',\mathbb{B}_{Y'}, \nu, \widetilde{T})$ forms a measure dynamical system. 

Let us define a map  $\eta$ 
from $X_{\sigma_0}$ to $Y'$. Take an element $z=x_1x_2\dots\in X_{\sigma_0}$. Then 
each prefix $x_1x_2\dots x_{\ell}$ with $\ell>3$ is a subword of 
the fix point $v$ of $\sigma_0$ starting with $a$. 
Therefore there is a word $y\in \{\lambda, a, ba, aba\}$ 
and $z_1\in X_{\sigma_0}$ such that $x=y_1\sigma_0(z_1)$. 
It is easy to see 
from (\ref{Delay}) that this $y$ and $z_1$ are unique. 

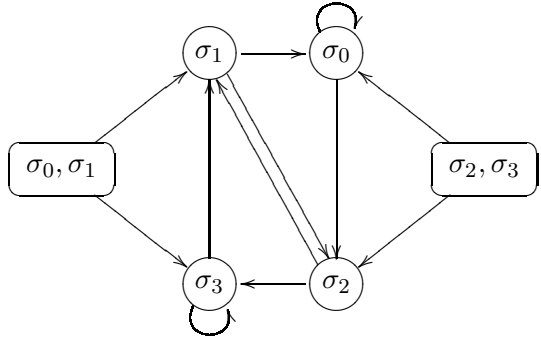
\begin{figure}
\[
\xymatrix{
  & *++[o][F]{\sigma_1} \ar@<0.5ex>[ddr] \ar[r] &  *++[o][F]{\sigma_0} \ar@(ul,ur) \ar[dd] & \\
*++[F-:<4pt>]{\sigma_0,\sigma_1}\ar[dr] \ar[ur] &   &  & *++[F-:<4pt>]{\sigma_2,\sigma_3} \ar[dl] \ar[ul]\\
  & *++[o][F]{\sigma_3} \ar@(dl,dr) \ar[uu] & *++[o][F]{\sigma_2} \ar[l] \ar@<0.5ex>[uul] & \\
}
\]
\caption{Forbidden suffix of $\mathbf{A}'$\label{Forbidden}}
\end{figure}

Iterating this we have $z_{i}=y_{i+1}\sigma(z_{i+1})$ with $y_i\in \{\lambda, a, ba, aba\}$,
$z_{i+1}\in X_{\sigma_0}$ and $z_{0}=z$. Thus we have for any $\ell$, 
\begin{eqnarray*}
z&=&y_1\sigma_0(y_2\sigma_0(y_3\sigma_0\dots y_{\ell}(\sigma_0(z_{\ell}))))\\
 &=&y_1 \sigma_0(y_2) \sigma_0^2(y_3) \dots \sigma_0^{\ell-1}(y_{\ell}) \sigma_0^{\ell}(z_{\ell}).
\end{eqnarray*}
Define a map from $\{\lambda,a,ba,aba\}$ to $\Z$ by
$$
\kappa(\lambda)=0, 
\kappa(a)=1,
\kappa(ba)=2,
\kappa(aba)=3.
$$
Then $z_{i}=y_{i+1}\sigma(z_{i+1})$ is equivalent to
$z_i=\sigma_{\kappa(y_{i+1})}(z_{i+1})$ and $z$ is represented as an $S$-adic limit:
$$
z=\lim_{\ell\rightarrow \infty}\sigma_{\kappa(y_1)} \circ 
\sigma_{\kappa(y_2)}\circ \dots \circ \sigma_{\kappa(y_{\ell})}(z_{\ell}).
$$
for $\ell=1,2,\dots$. 
This gives a {\bf multiplicative coding} $\mathbf{d'}:X_{\sigma_0} \rightarrow 
\{\sigma_0,\sigma_1,\sigma_2,\sigma_3\}^{\N}$. 
Let $\mathbf{A}'$ be the
points of $X_{\sigma_0}$ whose multiplicative coding does not end up 
in an infinite word produced by reading the vertex labels of 
Figure \ref{Forbidden}.
Let us associate to $z$ a
$2$-adic integer $\iota(z)=-\sum_{i=0}\kappa(y_i)2^{2i}\in \Z_2$. The map
$\iota$ is clearly bijective bi-continuous and the value $\iota(z)$ is
also called the multiplicative coding of $z$.
We write down first several iterates of $V$
on the fix point of $\sigma_0$, to illustrate the situation:
\begin{eqnarray*}
\sigma_0 \sigma_0 \sigma_0 \sigma_0 \dots &\stackrel{\iota}{\rightarrow}& -0000\dots\\
\sigma_3 \sigma_3 \sigma_3 \sigma_3 \dots &\stackrel{\iota}{\rightarrow}& -3333\dots\\
\sigma_2 \sigma_3 \sigma_3 \sigma_3 \dots &\stackrel{\iota}{\rightarrow}& -2333\dots\\
\sigma_1 \sigma_3 \sigma_3 \sigma_3 \dots &\stackrel{\iota}{\rightarrow}& -1333\dots\\
\sigma_0 \sigma_3 \sigma_3 \sigma_3 \dots &\stackrel{\iota}{\rightarrow}& -0333\dots\\
\sigma_3 \sigma_2 \sigma_3 \sigma_3 \dots &\stackrel{\iota}{\rightarrow}& -3233\dots
\end{eqnarray*}
One can see that the following commutative diagram (\ref{Additive0}) holds. 
\begin{equation}
\label{Additive0}
\begin{CD}
X_{\sigma_0} @>V>> X_{\sigma_0}\\
@V{\iota}VV              @V{\iota}VV\\
\Z_2  @>+1>> \Z_2
\end{CD}
\end{equation}
Therefore $(X_{\sigma_0},V)$ is topologically conjugate to the
$2$-adic odometer $(\Z_2,x\mapsto x+1)$. Here the consecutive digits $\{0,1\}$ 
in $\Z_2$ are glued together to give $\{0,1,2,3\}=\{0,1\}+2\{0,1\}$. 
Indeed, $\sigma_0$ satisfies the coincidence condition of height one
in the sense of Dekking \cite{Queffelec:87, Fogg:02} 
and above conjugacy is a consequence of this. 
$(\Z_2,x\mapsto x+1)$ is a translation of a compact group
$\Z_2$ which is minimal and uniquely ergodic with the Haar measure of $\Z_2$. 
Moreover one can confirm that $\iota$ preserves the measure and
$(X_{\sigma_0},V)$ and  $(\Z_2,x\mapsto x+1)$ are isomorphic through $\iota$
as measure dynamical systems. 
In view of (\ref{Delay}), 
we define 
$$
\xi(i)=\begin{cases} 0 & i=\lambda\\
                     2 & i=a\\
                     3 & i=ba\\
                     5 & i=aba
                     \end{cases}
$$
and the map $\eta:X_{\sigma_0} \rightarrow Y'$ by
\begin{equation}
\label{expansion}
\eta(x)=d_{\xi(y_1)}
+\frac {\zeta^{\xi(y_1)}}{\omega^2}\left(
d_{\xi(y_2)}+\frac {\zeta^{\xi(y_2)}}{\omega^2}\left(
d_{\xi(y_3)}+\frac {\zeta^{\xi(y_3)}}{\omega^2}
\dots. \right.\right.
\end{equation}
Then $\eta$ is clearly surjective, continuous, and measurable because
both $\tau$ and $\nu$ are Borel probability measures. 
Since the set of points with double addresses is on the open edge, 
the map $\eta$ is bijective from $\mathbf{A'}$ to $\mathbf{A}\cap T(\mathcal{Z})$.
Since $\widetilde{\mathbf{d}}(T(x))=V(\widetilde{\mathbf{d}}(x))$, we have 
a commutative diagram:
\begin{equation}
\label{Additive1}
\begin{CD}
\mathbf{A'} @>V>> \mathbf{A'}\\
@V{\eta}VV              @V{\eta}VV\\
 \mathbf{A}\cap T(\mathcal{Z})  @>\widetilde{T}>>  \mathbf{A}\cap T(\mathcal{Z}).
\end{CD}
\end{equation}
From Figure \ref{OpenEdge2}, it
is easy to see that
the set $\mathcal{P}$
of $\widetilde{T}$-periodic points in $Y'$ is measure zero by $\nu$, i.e.,
$\nu(\mathbf{A}\cap T(\mathcal{Z}))=\nu(Y'\cap T(\mathcal{Z}))=1$, because the
number of words of length $n$ in Figure \ref{OpenEdge2} is $O(2^n)$.
Similarly
as the Perron-Frobenius root of the substitution $\sigma_0$
is $4$ and the number of words of lengths $n$ in Figure \ref{Forbidden} 
are $O(2^n)$, we see that
$\tau(\mathbf{A'})=\tau(X_{\sigma_0})=1$. 
From (\ref{Additive1}) the pull back measure 
$\nu \circ \eta^{-1}$ of $X_{\sigma_0}$
is invariant by $V$, we have
$\tau=\nu \circ \eta^{-1}$ by unique ergodicity. 
Therefore by taking $\phi=\eta\circ \iota$, we have 
the commutative diagram (\ref{Additive}) with measure zero exceptions. 
Let $V'$ be a map from $X_{\sigma_0}$ to itself which acts
as the shift operator on the multiplicative coding $\mathbf{d'}$, i.e.,
$(\mathbf{d'}(V'(z))=\sigma_{n_2}\sigma_{n_3}\dots$ 
for $\mathbf{d'}(z)=\sigma_{n_1}\sigma_{n_2}\dots$.
Then
we see that
\begin{equation}
\label{Multiplicative1}
\begin{CD}
\mathbf{A'} @>V'>> \mathbf{A'}\\
@V{\eta}VV              @V{\eta}VV\\
 \mathbf{A}\cap T(\mathcal{Z})  @>S>>  \mathbf{A}\cap T(\mathcal{Z}).
\end{CD}
\end{equation}
and the commutative diagram (\ref{Multiplicative}) is valid but for measure zero
exceptions.
\qed

\begin{cor}
\label{UD}
Each aperiodic point $x\in \mathbf{A}\cap T(\mathcal{Z})$, 
the $\widetilde{T}$-orbit of $x$ is uniformly distributed in $Y'$
with respect to the self similar measure $\nu$. 
\end{cor}

\proof In the proof of Theorem \ref{Main} the map $\eta$ is
bijective form $\mathbf{A'}$ to $\mathbf{A}\cap T(\mathcal{Z})$.
Therefore if $x\in \mathbf{A}\cap T(\mathcal{Z})$, then there exists 
a unique element in $z\in X_{\sigma_0}$ with $\eta(z)=x$. Therefore there exist an
element $z_0\in \Z_2$ such that $\phi(z_0)=x$. The Haar measure $\mu_2$ on $Z_2$ 
is given by the values on the semi-algebra:
$$
\mu_2([c_0,c_1,\dots, c_{\ell-1}])=4^{-\ell}
$$
for each cylinder set
$[c_0,c_1,\dots, c_{\ell-1}]=\{ y\in \Z_2 \ |\ y\equiv \sum_{i=0}^{\ell-1} c_i 4^i \pmod{4^{\ell}}
\}$.
Since $(\Z_2, x\mapsto x+1)$ is uniquely ergodic, the assertion follows immediately
from the commutative diagram (\ref{Additive}).
\qed

Not all points in $Y'$ gives a dense orbit as we already 
mentioned that $\mathbf{A}\cap Y'$ is a proper dense subset of $Y'$. There
are many periodic points in $Y'$ as well.
This gives a good contrast to usual minimal topological dynamics given by 
a continuous map acting on a compact metrizable space. 

\begin{cor}
\label{UD2}
Each aperiodic point $x\in \mathbf{A}$, 
the $T$-orbit of $x$ is dense in the set $X$. 
\end{cor}

\proof
It is clear from the fact that 
$(T(\mathcal{Z}),\widetilde{T})$ is the induced system of $(L,T)$.
\qed
\bigskip

\begin{figure}
\begin{center}
\psfig{figure=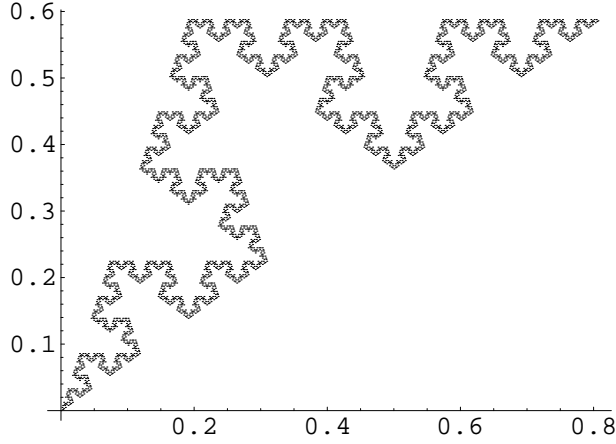}
\end{center}
\caption{The dual attractor $\mathcal{Y}$\label{Dual}}
\end{figure}

One can construct a dual expansion of the non-invertible dynamics $(Y',S)$
by the conjugate map $\phi:\zeta\rightarrow \zeta^2$
in ${\rm Gal} (\mathbb{Q}(\zeta)/\mathbb{Q})$ and then make a natural extension:
an invertible dynamics which contains $(Y',S)$.
The idea comes from
symbolic dynamics. We wish to construct the reverse expansion of (\ref{expansion})
to the other direction. To this matter, we compute in the following way:
$$
\frac{\omega^2(\eta(x)-d_{\xi(y_1)})}{\zeta^{\xi(y_1)}}=
d_{\xi(y_2)}+\frac{\zeta^{\xi(y_2)}}{\omega^2}\left(
d_{\xi(y_3)}+\frac{\zeta^{\xi(y_3)}}{\omega^2}\left(
d_{\xi(y_4)}+\frac{\zeta^{\xi(y_4)}}{\omega^2}\left(\dots \right.\right.\right.
$$
and 
$$
\frac{\omega^2}{\zeta^{\xi(y_2)}}\left(
\frac{\omega^2(\eta(x)-d_{\xi(y_1)})}{\zeta^{\xi(y_1)}}-d_{\xi(y_2)}\right)=
d_{\xi(y_3)}+\frac{\zeta^{\xi(y_3)}}{\omega^2}\left(
d_{\xi(y_4)}+\frac{\zeta^{\xi(y_4)}}{\omega^2}\left(\dots \right.\right.
$$
Therefore it is natural to introduce a left `expansion':
$$
\frac{\omega^2}{\zeta^{i_1}}\left(
\frac{\omega^2}{\zeta^{i_2}}\left(
\frac{\omega^2}{\zeta^{i_3}}\left(\left(\dots \right) -d_{i_3}\right)
-d_{i_2}\right) -d_{i_1}\right)
$$
with $i_k\in \{0,2,3,5\}$. As this expression does not converge,
we take the image of $\phi$ because $\phi(\omega)=-1/\omega$.
Let us denote by $u_{i_k}=\phi(d_{i_k})$. Then the expansion
$$
\frac {\zeta^{-2i_{1}}}{\omega^2}
\left(
\frac {\zeta^{-2i_{2}}}{\omega^2}
\left(
\frac {\zeta^{-2i_{3}}}{\omega^2} \left(\left( \dots \right) -u_{i_{3}}\right) \right) - u_{i_{2}}\right)- u_{i_{1}}
$$
converges and the closure of the set of such expansions gives 
a compact set $\mathcal{Y}$ depicted in figure \ref{Dual}. 

Of course
the set is an attractor of the IFS:
$$
\mathcal{Y}=
\frac {1}{\omega^2}(\mathcal{Y}-u_0) \cup
\frac {\zeta}{\omega^2}(\mathcal{Y}-u_2) \cup
\frac {\zeta^{-1}}{\omega^2}(\mathcal{Y}-u_3) \cup
\frac {1}{\omega^2}(\mathcal{Y}-u_5).
$$

Combining $\mathcal{Y}$ we can construct a natural extension of $(Y',S)$ as:
$$
Y'\times \mathcal{Y} \ni \left(\eta,\theta\right) \stackrel{\hat{S}}{\mapsto}
\left( \frac{(\eta-d_i)\omega^2}{\zeta^i}, 
\frac{\zeta^{-2i}(\theta-\phi(d_i))}{\omega^2}
\right) \in Y'\times \mathcal{Y}
$$
On the other hand $(\Z_2,\rho)$ have a natural extension: 
$$
\Z_2 \times [0,1) \ni (x,y)\stackrel{\hat{\rho}}{\mapsto} 
\left(\frac {x-(x \bmod{4})}4, \frac{y+(x \bmod{4})}4
\right) \in \Z_2\times [0,1)
$$
and two systems are isomorphic both as topological and measure theoretical dynamics:
\begin{equation}
\label{NaturalExt}
\begin{CD}
\Z_2 \times [0,1) @>\hat{\rho}>> \Z_2 \times [0,1)\\
@V\phi \times \phi'VV              @V\phi\times \phi'VV\\
Y'\times \mathcal{Y}  @>\hat{S}>>  Y'\times \mathcal{Y}.
\end{CD}
\end{equation}
where $\phi'$ is given as:
$$
\sum_{i=1}^{\infty} x_i 4^{-i} \mapsto g_{x_1}(g_{x_2}(g_{x_3}(\dots)))
$$
where $g_0(x)=(x-u_0)/\omega^2, g_1(x)=(x-u_2)\zeta/\omega^2,
g_2(x)=(x-u_3)\zeta^{-1}/\omega^2$ and $g_3(x)=(x-u_5)/\omega^2$.

From this `algebraic' natural extension construction, 
we can characterize purely $S$-periodic points in $Y'\cap \Q(\zeta)$. 

\begin{thm}
\label{PureExpansion}
A point $y$ in $Y'\cap \Q(\zeta)$
has purely periodic multiplicative coding with four digits $\sigma_0,
\sigma_2,\sigma_3,\sigma_5$ if and only if 
$(y,\phi(y))\in Y'\times \mathcal{Y}$.
\end{thm}

This is an analogy of the results \cite{Ito-Rao:05} for 
$\beta$-expansion. The proof below is on the same line.

\proof
As $\omega$ is an algebraic unit and $d_i\in \Z[\zeta]$, the denominator of $g_i(y)$
is the same as that of $y$ for $i=0,1,2,3$. Therefore the module $y\in \mathcal{M}=
\frac 1M\Z[\zeta]$
is stable by $g_i$ for some positive integer $M$. 
Note that points $y\in \mathcal{M}$
with $(y,\phi(y))\in  Y'\times \mathcal{Y}$ is finite, because 
$y,\phi(y)$ and their complex conjugates are bounded in $\C$. One can confirm that 
the map $\hat{S}$ becomes surjective from
$\mathcal{M}$ to itself. 
For a finite set, surjectivity 
implies bijectivity. Therefore a point $y\in \mathcal{M}$ with $(y,\phi(y))\in Y'\times
\mathcal{Y}$ produces a purely periodic orbit. 
On the other hand if $x$ has purely periodic multiplicative coding, it is
easy to see $(y,\phi(y))\in Y'\times \mathcal{Y}$. 
\qed

\section{Other self-similar systems} 
\label{sec:other_self_similar_systems}

Pisot scaling constants appear in several important dynamics.
For irrational rotations (2IET), 
it is well known that scaling constants of
self-inducing systems must be quadratic Pisot units. 
A typical example Figure \ref{fig:Int_exchange} was shown 
in the introduction.
They are computed 
by the continued fraction algorithm as fundamental units of quadratic number fields.
Poggiaspalla-Lowenstein-Vivald \cite{
PoggiaspallaLowensteinVivaldi:08} showed that 
the scaling constant must be an algebraic unit for self-inducing uniquely ergodic 
IET. When the scaling constant of self-inducing IET is a cubic Pisot unit,
we have further nice properties 
\cite{Arnoux-Rauzy:91,LowensteinPoggiaspallaVivaldi:07,LowensteinVivaldi:08}.

A necessary condition that $1$-dimensional
substitutive point sets give point diffraction is that the scaling constant 
is a Pisot number \cite{Bombieri-Taylor:87}. Suspension tiling dynamics
of such substitution is conjectured to have
 pure discrete spectrum if the characteristic polynomial of its
 substitution matrix is irreducible. 
For higher dimensional tiling dynamics the Pisot (or Pisot family) property is 
essential to have relatively dense point spectra, see for e.g. \cite{Solomyak:97, Lee-Solomyak:08}.

Pisot scaling properties seem to extend to the case of piecewise isometries. To conclude we present some examples, though we do not make a systematic study.

It is already observed in 
\cite{Kouptsov-Lowenstein-Vivaldi:02, 
Akiyama-Brunotte-Pethoe-Steiner:07} that 
Pisot scaling constants appear in our problem 
if $\theta$ is the $n$-th root of unity for $n=4,6,8,10,12$ in the same way as we did in $n=5$
but in a more involved manner. In each case they are quadratic Pisot units. 
What about if
$\lambda=-2\cos(\theta)$ is cubic? 
In this case, the dynamics of 
Conjecture \ref{Rotation} are embedded
into the piecewise affine mapping acting on $(\R/\Z)^4$ which is harder to visualize. Instead let us consider 
 formal analogies of piecewise isometries generated by cubic 
 $n$-th fold rotation in the plane.  At the expense of losing connection to
Conjecture \ref{Rotation}, we find many Pisot unit scaling constants!
Being an algebraic unit is natural and may be explained from
invertibility of dynamics. However we have 
no idea why the Pisot numbers turn up or even how to formulate these phenomena as a suitable
conjecture.

\subsection{Seven-fold} 
\label{sub:seven_fold}

\begin{figure}[htbp]
	\centering
		\includegraphics{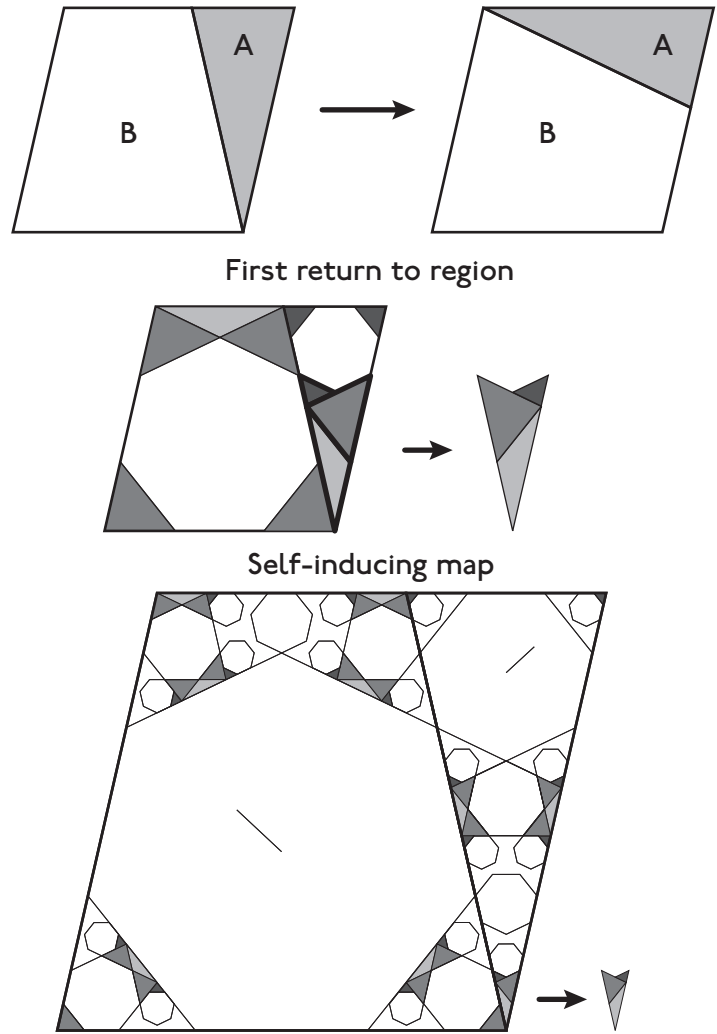}
	\caption{A seven-fold piecewise isometry.}
	\label{fig:final_7_fold_1}
\end{figure}

\begin{figure}[htbp]
	\centering
		\includegraphics{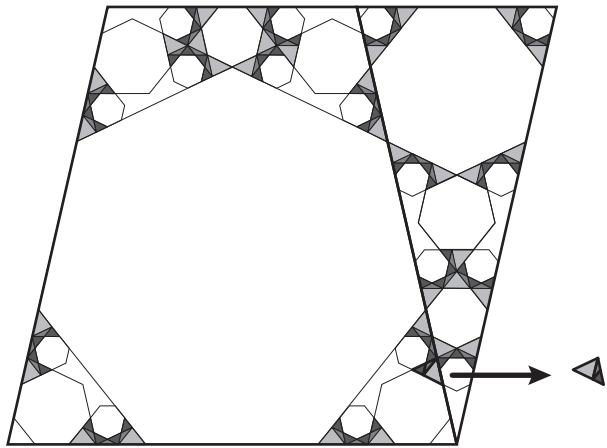}
	\caption{The regions remaining from the self-similarity shown in Figure~\ref{fig:final_7_fold_1}}
	\label{fig:final_7_fold_2}
\end{figure}

\begin{figure}[htbp]
	\centering
		\includegraphics{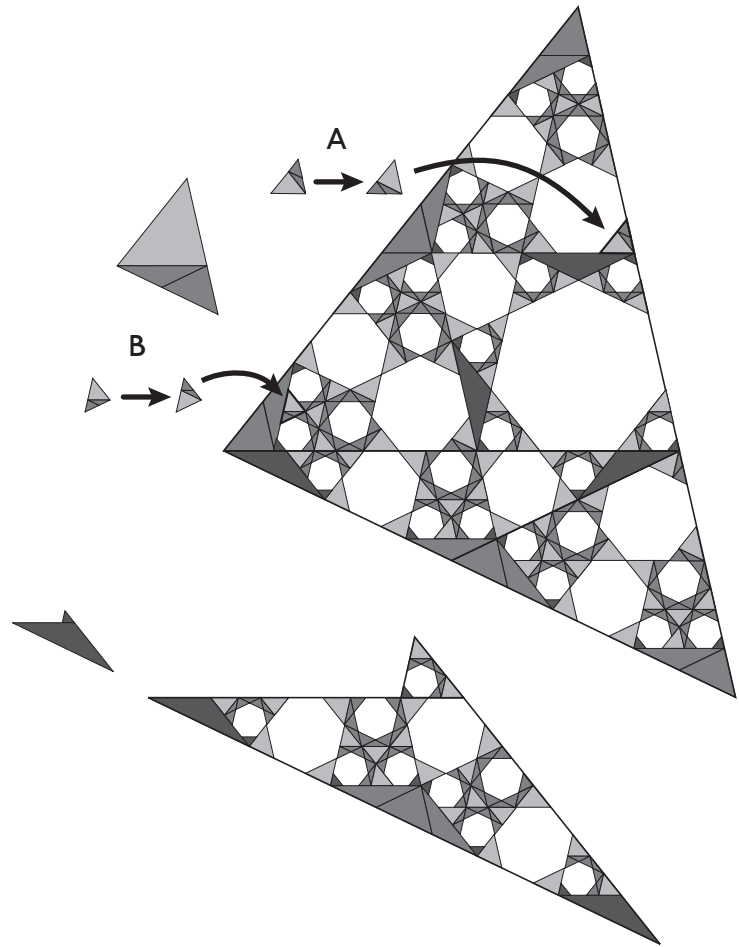}
	\caption{The substitution rule of the induced subsystem shown in Figure~\ref{fig:final_7_fold_2}.}
	\label{fig:final_7_fold_3}
\end{figure}

\begin{figure}[htbp]
	\centering
		\includegraphics{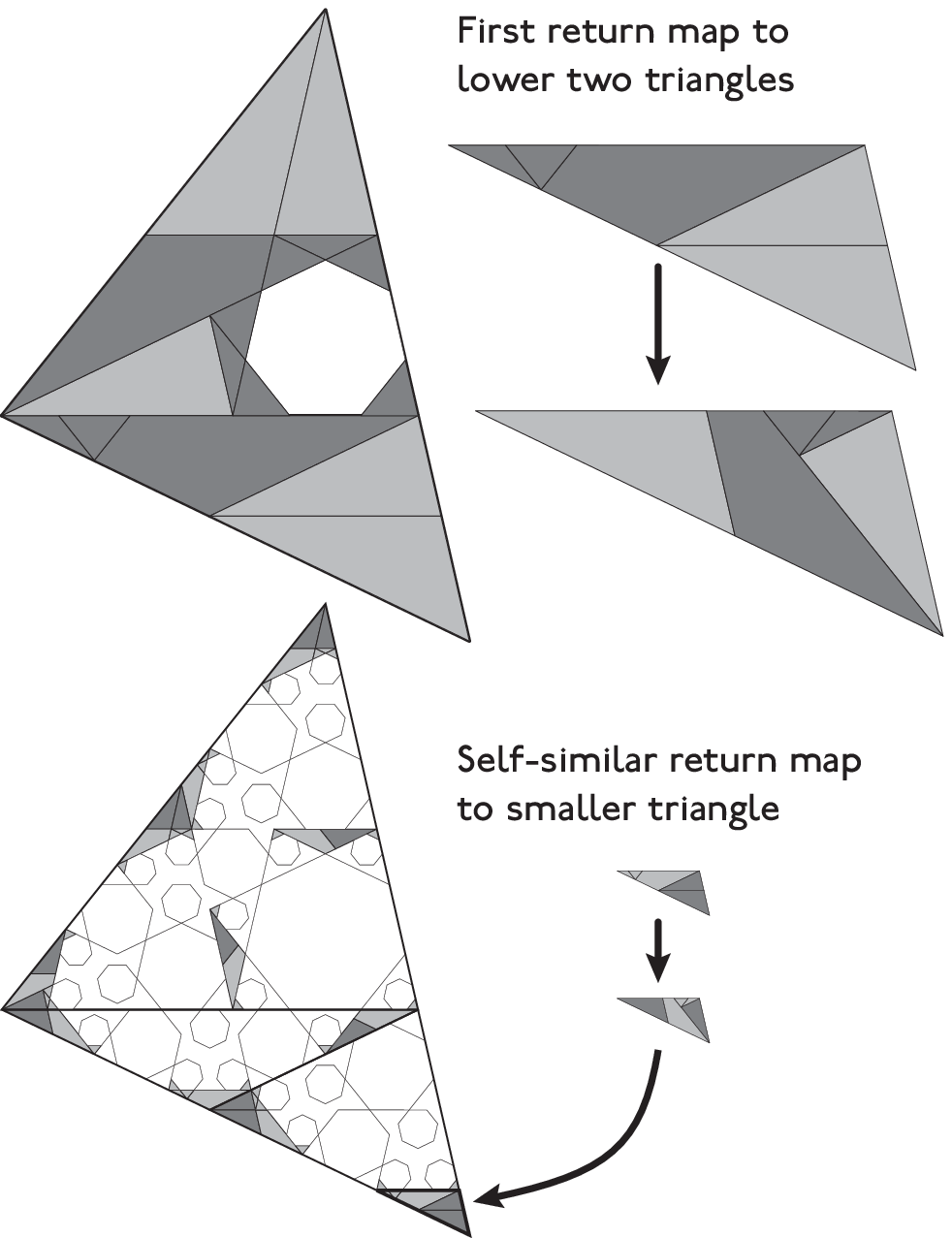}
	\caption{The final pieces of the structure of the piecewise isometry found in Figure~\ref{fig:final_7_fold_1}. }
	\label{fig:final_7_fold_4}
\end{figure}

We start with 7-fold case. 
Both pieces are rotated clockwise by $4 \pi/7$ as in Figure~\ref{fig:final_7_fold_1}. 
The triangle is rotated around A and the trapezium around B. 
The first return map to a region and a smaller region with the same first return map (up to scaling) are described. Unlike the five fold case, 
returning to the subregion does not cover the full region. 
A simple consequence is that there are infinitely many possible orbit closures for non-periodic orbits in the system. The scaling constant $\alpha \approx 5.04892$ 
is a Pisot number whose minimal polynomial is $x^3-6x^2+5x-1$.
Figure~\ref{fig:final_7_fold_2} shows how this remaining space can be filled in.
As this region is already a little small we will zoom in and now consider just this induced sub-system in Figure~\ref{fig:final_7_fold_3}. 
The smaller substitutions are easier to see as there are two scalings giving the same dynamics (A and B). The scaling constant $\beta \approx 16.3937$
for these subregions is the Pisot number associated to $x^3-17x^2+10x-1$.
The proof that the remaining substitutions work is shown in Figure~\ref{fig:final_7_fold_4}. 
The first return map to the two lower triangles is shown. The same dynamics occur on a smaller region. The orbit of the smaller region covers all the regions left out of Figure~\ref{fig:final_7_fold_3} and so the substitution rule from that figure is now complete.
The scaling constant for this triangle is $\alpha$. This gives an example of 
recursive tiling structure by Lowenstein-Kouptsov-Vivaldi 
\cite{Lowensein-Kouptsov-Vivaldi:04}.
Knowing that every aperiodic orbits are in one of the above self-inducing structures, we can show that

\begin{thm}
Almost all points of this $7$-fold lozenge have periodic orbits.
\end{thm}

\noindent
The argument is similar to that given around Figure \ref{Full}.
We easily find decreasing series $X_n$ of union of polygons satisfying
$\mu(\alpha X_{n+1})<\alpha^2 \mu(X_n)$ (or 
$\mu(\beta X_{n+1})<\beta^2 \mu(X_n)$)
which cover all self-inducing structures. 

The fundamental units of the maximal real subfield $\Q(\cos(4\pi/7))$ of 
the cyclotomic field $\Q(\zeta_7)$
are given by 
$b$ and $b-1$ where $b=1/(2\cos(3\pi/7))\approx 2.24698$. Here
$b$ is the Pisot number satisfying $x^3-2x^2-x+1$. 
We see that $\alpha=b^2$ and $\beta=b^4/(b-1)^2$ and thus
$\alpha$ and $\beta$ generates a subgroup of fundamental units
of $\Q(\cos(4\pi/7))$. Note that both $\sqrt{\alpha}=b$ and 
$\sqrt{\beta}=b^2/(b-1)$ are Pisot numbers but $b-1$ is not. 
Our piecewise isometry somehow selects Pisot units out of the unit group!

%
%
%

\subsection{Nine-fold} 
\label{sub:nine_fold}
\begin{figure}[htbp]
	\centering
		\includegraphics{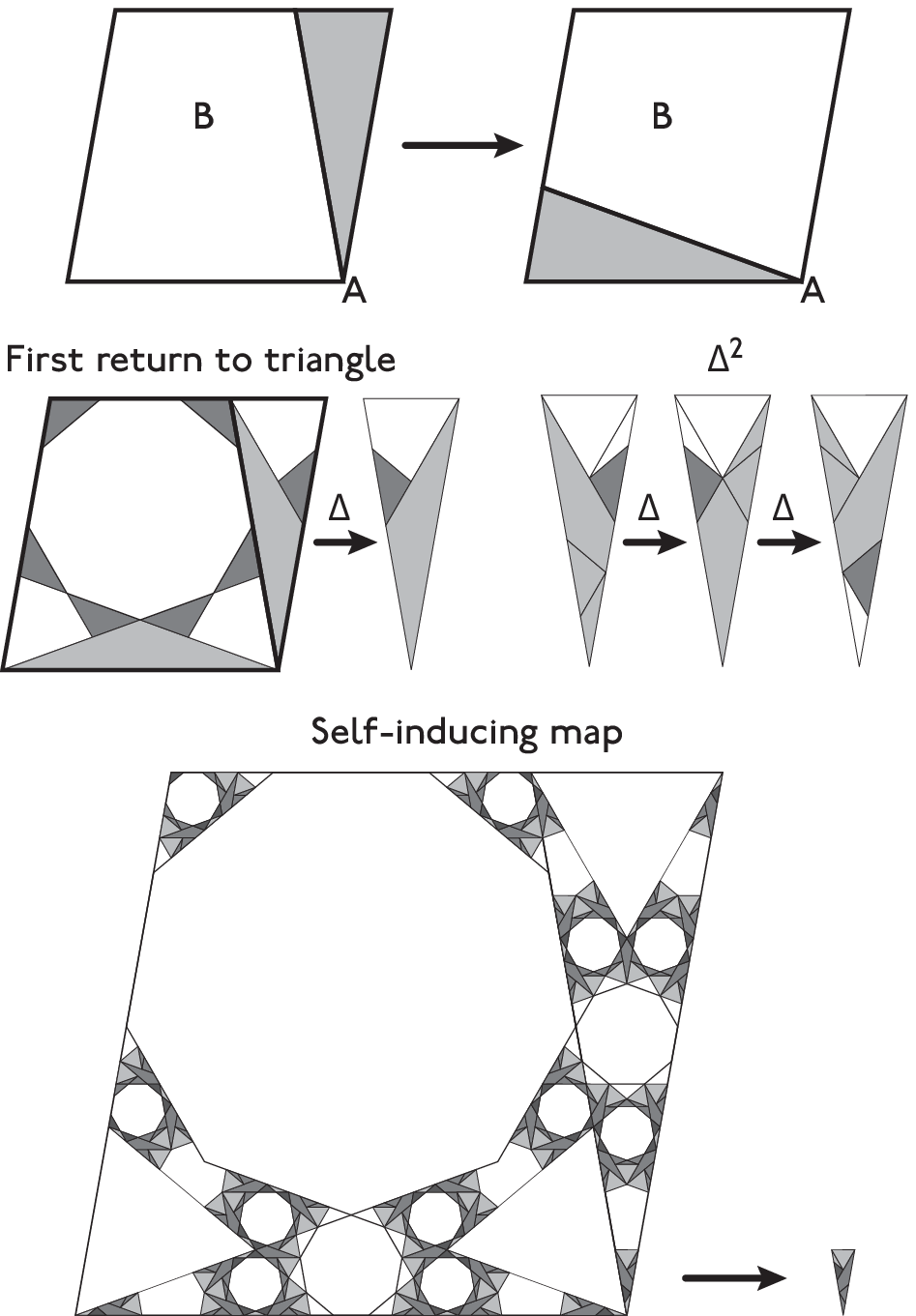}
	\caption{A nine fold piecewise isometry.}
	\label{fig:final_9_fold}
\end{figure}

The next example is 9-fold case in Figure \ref{fig:final_9_fold}. Both pieces are rotated anti-clockwise by $4\pi/9$, 
the triangle around A and the trapezium around B. The first return map ($\triangle$) to 
the triangle is also shown. In addition the same dynamics are found on a smaller 
piece of the map. 
Like the 7-fold shown in Figure~\ref{fig:final_7_fold_1} this does give a full description of the dynamics, but it is $\triangle^2$ not $\triangle$. 
The scaling constant $\gamma \approx 8.29086$ is a 
Pisot unit defined by $x^3-9x^2+6x-1$. 
Unfortunately in this case we were not able to find a complete description of the scaling structure.

 The fundamental units of $\Q(\cos(4\pi/9))$ are $b$ and $b^2-2b-1$ where
$b=1/(2\cos(4\pi/9))\approx 2.87939$ is 
a Pisot number given by $x^3-3x^2+1$. We have
$\gamma=b^2$ and are expecting to find another 
Pisot unit $b^2/(b^2-2b-1)\approx 5.41147$ (or its square) as
a scaling constant in this dynamics, which would give 
an analogy to the seven-fold case.

\bibliographystyle{amsplain}
\bibliography{./reflist}

\end{document}